\documentclass[a4paper,11pt,twoside]{article}
\usepackage{mfpaperstuff2}
\usepackage{tabu}
\usepackage[labelformat=simple]{caption}
\usepackage{ifthen}
\usepackage{chessfss}

\makeatletter
\newsavebox{\@brx}
\newcommand{\llparen}[1][]{\savebox{\@brx}{\(\m@th{#1(}\)}%
  \mathopen{\copy\@brx\kern-0.7\wd\@brx\usebox{\@brx}}}
\newcommand{\rrparen}[1][]{\savebox{\@brx}{\(\m@th{#1)}\)}%
  \mathclose{\copy\@brx\kern-0.7\wd\@brx\usebox{\@brx}}}
\makeatother

\begin{document}
\newcommand\pry[1]{\operatorname{pt}(#1)}
\newcommand\cbcx{canonical basis coefficient}
\newcommand\cbc{\cbcx\xspace}
\newcommand\cbcs{\cbcx s\xspace}
\newcommand\rdn[2]{\mathrm{D}_{#1#2}}
\newcommand\bp[2]{[#1:#2]}
\newcommand\hs{\tilde{\mathfrak{S}}_}
\newcommand\can[1]{\operatorname{G}(#1)}
\newlength\jmp\setlength\jmp{0.5\baselineskip}
\newcommand\Item[1][]{%
  \ifx\relax#1\relax  \item \else \item[#1] \fi
  \abovedisplayskip=0pt\abovedisplayshortskip=0pt~\vspace*{-\baselineskip}}
\newcommand\fg{\Yfillcolour{gray}}
\newcommand\fw{\Yfillcolour{white}}
\newcommand\halpha{\hat\alpha}
\newcommand\hbeta{\hat\beta}
\newcommand\hgamma{\hat\gamma}
\newcommand\isp\omega
\newcounter{casecount}\setcounter{casecount}0
\newcommand\nextcase{\addtocounter{casecount}1\arabic{casecount}.}
\newcommand\rfo[2]{\mbox{\medmuskip=0mu\thickmuskip=0mu$(#1\,|\,#2)$}}
\newlength\plx\setlength\plx{15pt}
\newcommand\bea[1]{\bar{\ttb}_{#1}}
\newcommand\hbea[1]{\ttb_{#1}}
\newcommand\link[3]{#1\stackrel i\leftrightarrow#2}
\newcommand\lslex{\ls_{\operatorname{lex}}}
\newcommand\gslex{\gs_{\operatorname{lex}}}
\newcommand\lscolex{\ls_{\operatorname{colex}}}
\newcommand\gscolex{\gs_{\operatorname{colex}}}
\newcommand\llex{<_{\operatorname{lex}}}
\newcommand\glex{>_{\operatorname{lex}}}
\newcommand\lcolex{<_{\operatorname{colex}}}
\newcommand\gcolex{>_{\operatorname{colex}}}
\newcommand\ope[3]{\llparen#1,#2\rrparen_{#3}}
\newcommand\lo[1]{\mathtt{a}_{#1}}
\newcommand\hi[1]{\mathtt{b}_{#1}}
\newcommand\lohi{\mathtt{c}}
\newcommand\dn[2]{\mathrm{d}_{#1#2}}
\newcommand\abd{abacus display\xspace}
\newcommand\abds{abacus displays\xspace}
\newcommand\bc[1]{$#1$-bar-core\xspace}
\newcommand\bcs[1]{$#1$-bar-cores\xspace}
\newcommand\baw{bar-weight\xspace}
\newcommand\bac{bar-core\xspace}
\newcommand\bacs{bar-cores\xspace}
\newcommand\bw[1]{$#1$-bar-weight\xspace}
\newcommand\hbc{\bc h}
\newcommand\hbcs{\bcs h}
\newcommand\hbw{\bw h}
\newcommand\pbw[2]{\calp_{#1,#2}}
\newcommand\pb[1]{\pbw{#1}2}
\newcommand\hsp{$h$-strict partition\xspace}
\newcommand\hsps{$h$-strict partitions\xspace}
\newcommand\hstr{\calp_h}
\newcommand\blkw[2]{\mathrm{B}_{#1,#2}}
\newcommand\blk[1]{\blkw{#1}2}
\newcommand\cbv{canonical basis vector\xspace}
\newcommand\cbvs{canonical basis vectors\xspace}
\newcommand\ba[1]{\lan#1\ran}
\newcommand\ppi[1]{\figsymbol{R}_{#1}}
\newcommand\nat[1]{\figsymbol{B}_{#1}}
\newcommand\shp[1]{\figsymbol{K}_{#1}}
\newcommand\flt[1]{\figsymbol{N}_{#1}}
\newcommand\xx[1]{\figsymbol{Q}_{#1}}
\newcommand\yy[1]{\figsymbol{p}_{#1}}
\newcommand\bb[1]{\Gamma(#1)}
\newcommand\ddd\partial
\newcommand\col{\operatorname{col}}
\newcommand\len[1]{\operatorname{len}(#1)}

\Yvcentermath0
\Yboxdim{12pt}

\title{Defect $2$ spin blocks of symmetric groups\\and \cbcs}
\runninghead{Defect $2$ spin blocks and \cbcs}

\toptitle


\begin{abstract}
This paper addresses the decomposition number problem for spin representations of symmetric groups in odd characteristic. Our main aim is to find a combinatorial formula for decomposition numbers in blocks of defect $2$, analogous to Richards's formula for defect $2$ blocks of symmetric groups.

By developing a suitable analogue of the combinatorics used by Richards, we find a formula for the corresponding ``$q$-decomposition numbers'', i.e.\ the \cbcs in the level\nobreakdash-$1$ $q$-deformed Fock space of type $A^{(2)}_{2n}$; a special case of a conjecture of Leclerc and Thibon asserts that these coefficients yield the spin decomposition numbers in characteristic $2n+1$. Along the way, we prove some general results on $q$-decomposition numbers. This paper represents the first substantial progress on canonical bases in type $A^{(2)}_{2n}$.
\end{abstract}

\tableofcontents

\linespread{1.016}\selectfont

\section{Introduction}

The most significant outstanding problem in the representation theory of the symmetric group $\sss m$ is the determination of the \emph{decomposition numbers}, describing the composition factors of the reduction of an ordinary irreducible representation modulo a prime. A complete solution to this problem seems to be far out of reach at present, but a wide variety of results are known dealing with special cases. One of these is Richards's combinatorial formula \cite[Theorem~4.4]{rich} giving the decomposition numbers for all blocks of symmetric groups of defect~$2$.

The theory of decomposition numbers for \emph{projective} representations of $\sss m$ (or equivalently, representations of a Schur cover $\hs m$) is much less advanced. The faithful representations of $\hs m$, i.e.\ those which do not descend to representations of $\sss m$, are called \emph{spin} representations of $\hs n$. Although the ordinary irreducible spin characters were classified by Schur in 1911 \cite{schu}, the corresponding representations were not constructed until 1990, by Nazarov \cite{naz}. As with the case of representations of $\sss n$, the combinatorics of partitions plays a central role in the theory. The modular theory of spin representations was initiated in the 1960s by Morris, who conjectured the block structure for spin representations; this conjecture was proved by Humphreys \cite{hump}. But a suitable parameterisation of the irreducible modular spin representations was not found until 2002, by Brundan and Kleshchev, who proved analogues of Kleshchev's ``modular branching rules'' describing the effect of inducing or restricting a simple module. Decomposition numbers for spin representations have been computed in degree at most $18$ \cite{moya2,bmo,maas}, but very few general results are known. These general results include the Brundan--Kleshchev regularisation theorem \cite{bkreg} and M\"uller's determination \cite{muller} of the decomposition numbers for blocks of defect $1$. In this paper we address spin blocks of defect $2$, in the hope of finding a spin version of Richards's formula.

In fact most of this paper is concerned with quantum algebra. For every Kac--Moody algebra~$\fkg$ of classical affine type, Kashiwara, Miwa, Petersen and Yung \cite{kmpy} construct a \emph{$q$-deformed Fock space} of level~$1$; this is a module for the quantum group $U_q(\fkg)$. In the case where $\fkg$ is of type $A_{2n}^{(2)}$, Leclerc and Thibon \cite{lt} studied this Fock space further, introducing partition combinatorics and drawing a connection with spin representations of $\hs m$ in characteristic $h=2n+1$; this connection revolves around the fact that the action of the standard generators of $U_q(\fkg)$ corresponds to Morris's branching rules describing induction and restriction of spin representations between $\hs m$ and $\hs{m+1}$. The submodule of the Fock space generated by the empty partition is called the \emph{basic representation}, and possesses an important basis called the \emph{canonical basis}. The coefficients expressing canonical basis elements in terms of the standard basis for the Fock space are called ``$q$-decomposition numbers'', in view of a conjecture by Leclerc and Thibon \cite[Conjecture~6.2]{lt} that (after specialising at $q=1$ and suitable rescaling) these coefficients coincide with decomposition numbers for spin representations of $\hs m$ in characteristic $h$, provided $h$ is sufficiently large.

The main results of the present paper (\cref{mainwt1,mainwt2}) are combinatorial formul\ae{} for the $q$-decomposition numbers corresponding to spin blocks of $\hs n$ of defect $1$ or $2$; along the way, we prove a parity theorem, describing when a $q$-decomposition number is an odd or an even function of $q$. Our formula for defect~$1$ (combined with M\"uller's results) shows that the Leclerc--Thibon conjecture holds for blocks of defect~$1$. Our formula for defect $2$ is very similar in spirit to Richards's formula, though there is some ``exceptional'' behaviour for up to three canonical basis vectors in each block. In fact, the Leclerc--Thibon conjecture in its original formulation is now known to be false, as it predicts negative decomposition numbers (the author is grateful to Shunsuke Tsuchioka for providing this information). However, we expect that it is true for blocks of defect $2$, so that our formula specialises to give a formula for the decomposition numbers for spin blocks of defect $2$.

The results in this paper appear to be the first significant results on $q$-decomposition numbers in type $A^{(2)}_{2n}$. It is to be expected that similar results will hold in other affine types, and can be proved via the same techniques. After proving our main results, we provide a brief discussion of corresponding results for type $A^{(2)}_{2n+1}$, and their relationship to the results for type $A^{(2)}_{2n}$, which partly explains the exceptional behaviour seen in some of the canonical basis vectors.

Our main technique is to exploit the combinatorics of partitions, reconciling the action of the quantum group on the Fock space with the combinatorial notions underlying the formula, in particular leg lengths and the dominance order. We build on the work of Kessar and Schaps \cite{ke,kesc} to derive properties of \emph{Scopes--Kessar pairs} of blocks to enable inductive proofs of our main results. The techniques we develop can be applied to blocks of higher defect, and modified to provide results for other Kac--Moody types.

\begin{ack}
The research in this paper would not have been possible without extensive calculations using GAP \cite{gap}.
\end{ack}

\section{Combinatorial background}\label{combsec}

\begin{basicbox}
Throughout this paper $h$ denotes an odd integer greater than $1$, and we write $h=2n+1$.
\end{basicbox}

In this section we outline the combinatorial set-up underlying both the Fock space of type $A^{(2)}_{h-1}$ and the modular spin representations of symmetric groups.

\subsection{\hsps, cores and blocks}\label{hspsec}

A \emph{partition} is an infinite weakly decreasing sequence $\la=(\la_1,\la_2,\dots)$ of non-negative integers with finite sum. We write $\card\la=\la_1+\la_2+\dots$, and say that $\la$ is a partition of $|\la|$. The integers $\la_1,\la_2,\dots$ are called the \emph{parts} of $\la$, and the number of positive parts of $\la$ is called the \emph{length} of $\la$, written $\len\la$. We may write $a\in\la$ or say that $\la$ \emph{contains} $a$ if $\la_r=a$ for some $r$. When writing partitions, we usually group together equal parts with a superscript and omit the trailing zeroes, and we write the unique partition of $0$ as $\varnothing$. We say that $\la$ is \emph{strict} if $\la_r>\la_{r+1}$ for all $r<\len\la$.

The \emph{Young diagram} of a partition $\la$ is the set
\[
[\la]=\lset{(r,c)\in\bbn^2}{c\ls\la_r}
\]
whose elements we call the \emph{nodes} of $\la$. In general, a node means an element of $\bbn^2$. We draw Young diagrams as arrays of boxes using the English convention, in which $r$ increases down the page and~$c$ increases from left to right. The \emph{conjugate partition} to $\la$ is the partition $\la'$ whose Young diagram is obtained by reflecting $[\la]$ on the main diagonal; that is, $\la'_r=\max\lset{c\gs1}{\la_c\gs r}$.

The \emph{dominance order} is a partial order $\domby$ defined on the set of partitions of a given size by
\[
\la\domby\mu\qquad\Longleftrightarrow\qquad\la_1+\dots+\la_r\ls\mu_1+\dots+\mu_r\text{ for all }r\gs0.
\]
We will use a well-known alternative characterisation of the dominance order in term of conjugate partitions.

\begin{lemma}[\xcite{jk}{Lemma 1.4.11}]\label{conjdom}
If $\la$ and $\mu$ are partitions, then $\la\domby\mu$ \iff $\la'\dom\mu'$.
\end{lemma}

We will also need two total orders on partitions. Given partitions $\la,\mu$, we write $\la\llex\mu$ if there is $r$ such that $\la_r<\mu_r$ while $\la_s=\mu_s$ for all $s<r$. We say that $\la\lcolex\mu$ if there is $r$ such that $\la_r>\mu_r$ while $\la_s=\mu_s$ for all $s>r$. Then $\lslex$ and $\lscolex$ are total orders (called the \emph{lexicographic} and \emph{colexicographic} orders) which both refine the dominance order on partitions of a given size.

Finally we introduce some natural set-theoretic notation for partitions. Suppose $\la$ and $\mu$ are partitions.
\begin{itemize}
\item
We write $\la\sqcup\mu$ for the partition obtained by combining the parts of $\la$ and $\mu$ and arranging them into decreasing order.
\item
We write $\la\cap\mu$ for partition in which the number of parts equal to $a$ is the smaller of the number of parts of $\la$ equal to $a$ and the number of parts of $\mu$ equal to $a$, for each $a$.
\item
If $\mu$ is strict and $\mu_r\in\la$ for each $r$, we define $\la\setminus\mu$ to be the partition obtained by deleting one copy of $\mu_r$ from $\la$, for each $r$.
\end{itemize}

\medskip
Now we introduce the odd integer $h$ into the combinatorics. A partition $\la$ is \emph{$h$-strict} if for every $r$ either $\la_r>\la_{r+1}$ or $\la_r\equiv0\ppmod h$. An \hsp is \emph{restricted} if for every $r$ either $\la_{r+1}>\la_r-h$ or $\la_{r+1}=\la_r-h\nequiv0\ppmod h$. Throughout this paper we write $\hstr$ for the set of all \hsps.

For example, the $3$-strict partitions of $8$ are
\[
(8),(7,1),(6,2),(5,3),(5,2,1),(4,3,1),(3^2,2),
\]
and of these only the last three are restricted.

The \emph{residue} of a node $(r,c)$ is the smaller of the residues of $c-1$ and $-c$ modulo $h$. A node of residue $i\in\{0,1,\dots,n\}$ is called an \emph{$i$-node}.

For example, if we take $h=5$ and $\la=(11,8,6,5^2)$, the residues of the nodes of $\la$ are given in the following diagram.
\[
\young(01210012100,01210012,012100,01210,01210)
\]

Given $\la\in\hstr$ and $0\ls i\ls n$, let $\mu$ be the smallest \hsp such that $[\mu]\subseteq[\la]$ and $[\la]\setminus[\mu]$ consists entirely of $i$-nodes. These nodes are called the \emph{removable $i$-nodes} of $\la$. Note that in the case $i=0$, a removable $i$-node of $\la$ might not be a removable node in the conventional sense. For example, referring to the diagram above, the removable $0$-nodes of $(11,8,6,5^2)$ are $(1,10)$, $(1,11)$, $(3,6)$ and $(5,5)$ when $h=5$.

Similarly, the \emph{addable $i$-nodes} of $\la$ are the $i$-nodes that can be added to $\la$ (possibly together with other $i$-nodes) to create a larger \hsp. For example, when $h=5$ the addable $1$-nodes of the partition $\la=(11,8,6,5^2)$ above are $(1,12)$, $(2,9)$ and $(3,7)$.

The \emph{$h$-content} of a partition $\la$ is the multiset of residues of the nodes of $\la$. For example, when $h=5$ the partition $(11,8,6,5^2)$ has fifteen $0$-nodes, thirteen $1$-nodes and seven $2$-nodes, so we write its $5$-content as $\{0^{15},1^{13},2^7\}$.

Now we introduce \hbcs and blocks. Suppose $\la\in\hstr$. \emph{Removing an $h$-bar} from $\la$ means constructing a smaller \hsp by doing one of two things:
\begin{itemize}
\item
replacing a part $\la_r\gs h$ with $\la_r-h$ and reordering the parts into decreasing order;
\item
removing two parts which sum to $h$.
\end{itemize}
An \hsp $\la$ is called an \emph{\hbc} if it is not possible to remove an $h$-bar from $\la$. In general, the \emph{\hbc} of $\la\in\hstr$ is the \hbc obtained by repeatedly removing $h$-bars until it is not possible to remove any more. It is an easy exercise to show that the \hbc of $\la$ is well-defined, and the \emph{\hbw} is the number of $h$-bars removed to reach the \hbc.

Suppose $\tau$ is an \hbc and $w\gs0$. We define $\pbw\tau w$ to be the set of all \hsps with \hbc $\tau$ and \hbw $w$.

For example, suppose $\la=(9,6,3,1)$. Then the \bc5 of $\la$ is $(3,1)$, and its \bw5 is $3$, as we see from the following diagrams.
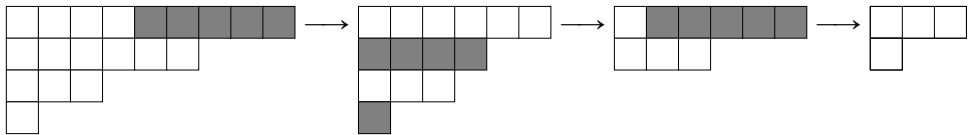
\[
\gyoung(;;;;!\fg;;;;;,!\fw;;;;;;,;;;,;)\longrightarrow\gyoung(;;;;;;,!\fg;;;;,!\fw;;;,!\fg;)\longrightarrow\gyoung(!\fw;!\fg;;;;;,!\fw;;;)\longrightarrow\yng(3,1)
\]

Later we will need the following results.

\begin{propn}[\xcite{moya}{Theorem 5}]\label{contcore}
Suppose $\la,\mu\in\hstr$ with $\card\la=\card\mu$. Then $\la$ and $\mu$ have the same \hbc \iff they have the same $h$-content.
\end{propn}

\begin{lemma}\label{non2core}
Suppose $\tau$ is an \hbc. Then one of the following occurs.
\begin{itemize}
\item
$\tau=(l,l-1,\dots,1)$ for some $0\ls l\ls n$.
\item
$\tau$ has a removable $i$-node, for some $i\neq0$.
\item
$\tau$ has at least two removable $0$-nodes.
\end{itemize}
\end{lemma}

\begin{pf}
Let $l=\len\la$. If $\tau=(l,l-1,\dots,1)$, then necessarily $l\ls n$, since otherwise $\tau$ would include the parts $n$ and $n+1$, so would not be an \hbc. So assume that $\tau\neq(l,l-1,\dots,1)$. Assume also that $\tau$ has no removable nodes of any non-zero residue. By assumption $l>0$, so the node $(l,\tau_l)$ is removable, and is therefore a removable $0$-node. This is particular means that $\tau_l=1$ (since otherwise it would be possible to remove an $h$-bar from $\tau$), and now the assumption $\tau\neq(l,l-1,\dots,1)$ means that there is $r<l$ such that $\tau_r-\tau_{r+1}\gs2$. So the node $(r,\tau_r)$ is also a removable $0$-node. So $\tau$ has at least two removable $0$-nodes.
\end{pf}

We end this section with some more notation which we shall use repeatedly: if $\tau$ is a strict partition and $x<y$ are integers, we write $\ope xy\tau$ for the number of parts of $\tau$ lying strictly between $x$ and~$y$. As a special case of this, we set $\bb\tau=\ope0h\tau$.

\subsection{The abacus}\label{absec}

Abacus notation for partitions was introduced by James, and has proved to be a valuable tool in the combinatorial modular representation theory of symmetric groups. A different abacus notation for strict partitions was introduced by Bessenrodt, Morris and Olsson to play the analogous role in the theory of spin representations. Here we introduce an alternative abacus notation which works better for our purposes (this is also used by Yates in \cite{deanyates}).

We take an abacus with $h$ vertical runners numbered $-n,1-n,\dots,n-1,n$ from left to right. On runner $i$ we mark positions labelled with the non-zero integers in $i+h\bbz$ increasing down the runner, so that (if $a\nequiv n\ppmod h$) position $a+1$ appears directly to the right of position $a$. For example, if $h=9$, the abacus is drawn as follows.
{\small
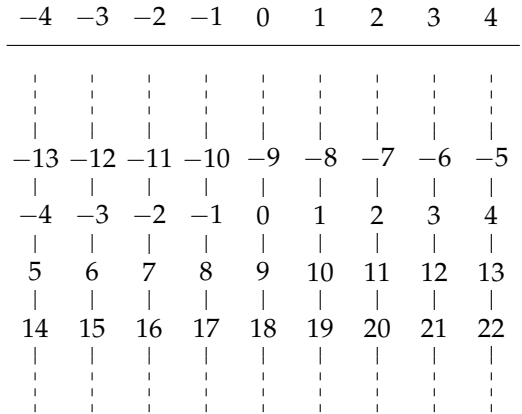
\[
\begin{tikzpicture}[scale=.75]
\foreach\x in{-4,...,4}{\draw(\x,3.5)node{$\x$};\draw(\x,-2.5)--++(0,4);\draw[dashed](\x,-2.5)--++(0,-1);\draw[dashed](\x,1.5)--++(0,1);}
\foreach\x in{1,...,4}{\draw(\x,0)node[fill=white]{$\x$};\draw(-\x,0)node[fill=white]{$-\x$};};
\draw(-4.5,3)--++(9,0);
\foreach\x in{-13,...,-5}\draw(\x+9,1)node[fill=white]{$\x$};
\foreach\x in{5,...,13}\draw(\x-9,-1)node[fill=white]{$\x$};
\foreach\x in{14,...,22}\draw(\x-18,-2)node[fill=white]{$\x$};
\end{tikzpicture}
\]
}
Now given a strict partition $\la$ of length $l$, the \emph{\abd} for $\la$ is obtained by placing beads in positions $\la_1,\dots,\la_l$ and in all negative positions except $-\la_1,\dots,-\la_l$. We place a $\times$ where position~$0$ would be.

For example, if $h=9$ and $\la=(15,11,5,4,2,1)$, the \abd for $\la$ is as follows.
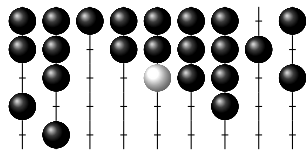
\[
\abacus(bbbbbbbnb,bbnbbbbbn,nbnnxbbnb,bnnnnnbnn,nbnnnnnnn)
\]

Here we have used a convention we shall apply throughout the paper: whenever we show an \abd (or just a portion of an \abd consisting of certain chosen runners), all positions above those shown are understood to be occupied, and all positions below those shown are understood to be unoccupied.

We shall also occasionally consider the abacus display of a partition which is $h$-strict but not strict. In this case, if a part $ah$ occurs $t$ times in $\la$, we regard the abacus as having $t$ beads at position $ah$, and $t$ empty spaces at position $-ah$; we depict this by labelling the bead at position $ah$ and the empty space at position $-ah$ with the integer $t$.

The effect of adding a node to an \hsp is easy to see on the abacus. The following \lcnamecref{addnodesabacus} follows from the definitions.

\begin{lemma}\label{addnodesabacus}
Suppose $\la,\mu\in\hstr$, and that $\mu$ is obtained from $\la$ by adding an $i$-node.
\begin{enumerate}
\item
If $1\ls i\ls n$, then the \abd for $\mu$ is obtained from the \abd for $\la$ by moving a bead from position $b$ on runner $i$ to runner position $b+1$, and simultaneously moving a bead from position $-b-1$ to position $-b$.
\item
If $i=0$, then the \abd for $\mu$ is obtained from the \abd for $\la$ either by moving a bead from position $b$ on runner $0$ to position $b+1$ and simultaneously moving a bead from position $-b-1$ to position $-b$, or moving a bead from position $-1$ to position $1$.
\end{enumerate}
\end{lemma}

The \abd also makes it easy to visualise removal of $h$-bars and construction of the \hbc of an \hsp. Suppose $\la$ is an \hsp from which we can remove an $h$-bar. There are three ways we do this, and we consider the effect on the \abd in each case.
\begin{itemize}
\item
We can replace a part $a+h$ with $a$, where $a\gs1$. In this case on the abacus we move the bead at position $a+h$ to position $a$, and we move the bead at position $-a$ to position $-a-h$.
\item
We can delete two parts $a$ and $h-a$, where $1\ls a<h$. In this case we move the beads at positions $a,h-a$ to positions $a-h,-a$.
\item
We can delete the part $h$. In this case we move the bead at position $h$ to position $-h$.
\end{itemize}

We see that in each case, removing an $h$-bar involves moving beads up their runners into unoccupied positions. As a consequence, we find that the \abd for the \hbc of $\la$ may be obtained by moving all beads up their runners as far as they will go. In particular, we have the following lemma.

\begin{lemma}\label{coreflush}
Suppose $\tau\in\hstr$. Then $\tau$ is an \hbc \iff every bead in the \abd for $\tau$ has a bead immediately above it.
\end{lemma}

Returning to the partition in the above example, we obtain the \abd for the $9$-bar-core $(11,6,2,1)$.
\[
\abacus(bbbbbbbbb,bbnbbbbnb,bbnnxbbnn,nbnnnnbnn,nnnnnnnnn)
\]

\section{The $q$-Fock space}\label{focksec}

Now we introduce the background we shall need from quantum algebra. This is essentially taken from the paper \cite{lt} by Leclerc and Thibon; note, however, that the residues used there are the opposite of ours: a node of residue $i$ in this paper has residue $n-i$ in \cite{lt}.

\subsection{The quantum algebra and the Fock space}

We consider the quantum group $\calu=U_q(A^{(2)}_{h-1})$ associated to the generalised Cartan matrix of type $A_{h-1}^{(2)}$. This comes with the usual Kac--Moody set-up of simple roots $\alpha_0,\dots,\alpha_n$ and fundamental weights $\La_0,\dots,\La_n$. We refer to the book by Hong and Kang \cite{hoka} for the necessary background on quantised Kac--Moody algebras.

$\calu$ has standard generators $e_i,f_i,t_i$ for $i=0,\dots,n$. We define
\[
q_i=
\begin{cases}
q&\text{if }i=0\\
q^2&\text{if }0<i<n\\
q^4&\text{if }i=n
\end{cases}
\]
and then set (for $k\gs0$)
\[
[k]_i=\frac{q_i^k-q_i^{-k}}{q_i-q_i^{-1}},\qquad[k]^!_i=[1]_i[2]_i\dots[k]_i,\qquad e_i^{(k)}=\frac{e_i^k}{[k]^!_i},\qquad f_i^{(k)}=\frac{f_i^k}{[k]^!_i}.
\]

The \emph{$q$-deformed Fock space} $\scrf$ (of level $1$) is a vector space over $\bbc(q)$ with $\hstr$ as a basis. This space is naturally a module for $\calu$, and we can give combinatorial rules for the action of the divided powers~$e_i^{(k)}$ and $f_i^{(k)}$. (Note that Leclerc and Thibon only give the actions of $e_i$ and $f_i$, and they express them in terms of straightening rules; but our rules can be easily deduced from the rules in \cite{lt}.)

First we consider $f_i^{(k)}$. Suppose $\la,\mu\in\hstr$ such that $[\mu]\supseteq[\la]$ and $[\mu]\setminus[\la]$ consists of $a$ nodes of residue $i$; then we write $\la\stackrel{a:i}\longrightarrow\mu$, and we define a coefficient $N_i(\la,\mu)$ as follows: let $s$ be the sum, over all nodes $(r,c)$ of $[\mu]\setminus[\la]$, of the number of addable $i$-nodes of $\mu$ to the left of $(r,c)$ minus the number of removable $i$-nodes of $\la$ to the left of $(r,c)$. Further, if $i=0$, let $M$ be the set of integers $m\gs1$ such that there is a node of $[\mu]\setminus[\la]$ in column $mh+1$ but not in column $mh$; for each $m\in M$, let~$b_m$ be the number of times $mh$ occurs as a part of $\la$, and set $f_m=1-(-q^2)^{b_m}$.

Now define
\[
N_i(\la,\mu)=\begin{cases}
q_i^s\prod_{m\in M}f_m&\text{if }i=0\\
q_i^s&\text{if }i\neq0.
\end{cases}
\]
Then
\[
f^{(a)}_i\la=\sum_{\la\stackrel{a:i}\longrightarrow\mu}N_i(\la,\mu)\mu.
\]

\begin{eg}
Take $h=5$ and $\la=(5,4)$. Let us calculate $f_0^{(k)}\la$ for each $k$. $\la$ has three addable $0$-nodes, namely $(1,6)$, $(2,5)$ and $(3,1)$. Applying the formula above, we obtain
\begin{align*}
f_0\la&=(5,4,1)+q(5^2)+(q^4+q^2)(6,4),\\
f_0^{(2)}\la&=(5^2,1)+(q^3+q)(6,4,1)+q^2(6,5),\\
f_0^{(3)}\la&=(6,4,1),\\
f_0^{(k)}\la&=0\quad\text{for all }k\gs4.
\end{align*}
\end{eg}

The rule for the action of $e_i$ is similar. We define a coefficient $N_i(\la,\mu)$ whenever $\la\stackrel{a:i}\longleftarrow\mu$ as follows: let $s$ be the sum, over all nodes $(r,c)$ of $[\la]\setminus[\mu]$, of the number of removable $i$-nodes of $\mu$ to the right of $(r,c)$ minus the number of addable $i$-nodes of $\la$ to the right of $(r,c)$. Further, if $i=0$, let $M$ be the set of integers $m\gs1$ such that there is a node of $[\la]\setminus[\mu]$ in column $mh$ but not in column $mh+1$; let~$b_m$ be the number of times $mh$ occurs as a part of $\la$, and set $f_m=1-(-q^2)^{b_m}$.

Now define
\[
N_i(\la,\mu)=\begin{cases}
q_i^s\prod_{m\in M}f_m&\text{if }i=0\\
q_i^s&\text{if }i\neq0.
\end{cases}
\]
Then
\[
e^{(a)}_i\la=\sum_{\la\stackrel{a:i}\longleftarrow\mu}N_i(\la,\mu)\mu.
\]

The Fock space $\scrf$ has a weight space decomposition $\scrf=\bigoplus_{\alpha}\scrf_{\La-\alpha}$, with $\alpha$ ranging over positive roots. If $\alpha$ is the root $\sum_{i=0}^na_i\alpha_i$ with each $a_i$ non-negative, then the weight space $\scrf_{\La_0-\alpha}$ is spanned by the \hsps $\la$ with $h$-content $\{0^{a_0},\dots,n^{a_n}\}$. So by \cref{contcore} two partitions lie in the same weight space \iff they have the same \bac and \baw. Given a \bac $\tau$ and an integer $w\gs0$, we write $\blkw\tau w$ for the weight space spanned by the \hsps with \bac $\tau$ and \baw $w$. We refer to $\blkw\tau w$ as the \emph{block} with \bac $\tau$ and \baw $w$; this helps us to avoid over-taxing the word ``weight'', and keeps in mind the connection with blocks in the sense of modular representation theory. The aim of this paper is to study blocks with small \baw.

\subsection{The canonical basis}

Now let $V_0$ denote the submodule of $\scrf$ generated by the empty partition $\varnothing$. This submodule is isomorphic to the irreducible highest-weight $\calu$-module $V(\La_0)$, and possesses a \emph{canonical basis}, defined as follows. The \emph{bar involution} is the $\bbc(q+q^{-1})$-linear involution $v\mapsto \widebar v$ on $V_0$ defined by $\widebar\varnothing=\varnothing$, $\widebar{e_iv}=e_i\widebar v$ and $\widebar{f_iv}=f_i\widebar v$. We say that a vector $v\in V_0$ is \emph{bar-invariant} if $v=\widebar v$; this means that $v$ can be written as a linear combination, with coefficients lying in $\bbc(q+q^{-1})$, of vectors of the form $f_{i_1}\dots f_{i_r}\varnothing$. For each restricted \hsp $\mu$, there is a unique vector $\can\mu\in\scrf$ with the following properties.
\begin{enumerate}[label=(CB\arabic*),leftmargin=*]
\item
$\can\mu$ is bar-invariant.
\item
When we write $\can\mu=\sum_\la\dn\la\mu\la$, the coefficient $\dn\mu\mu$ equals $1$ while all the other coefficients~$\dn\la\mu$ are polynomials divisible by $q$.
\end{enumerate}
The vector $\can\mu$ is called the \emph{canonical basis vector} corresponding to $\mu$, and the set
\[
\lset{\can\mu}{\mu\in\hstr}
\]
is the \emph{canonical basis} of $V_0$. Every bar-invariant vector in $V_0$ is a linear combination, with coefficients in $\bbc(q+q^{-1})$, of canonical basis vectors. The canonical basis can be computed recursively via the \emph{LT algorithm} \cite[Section~4]{lt}.

The coefficients $\dn\la\mu$ are called \emph{\cbcs} or \emph{$q$-decomposition numbers}, and they satisfy the following additional property:
\begin{enumerate}[resume,label=(CB\arabic*),leftmargin=*]
\item
If $\dn\la\mu\neq0$, then $\mu\domby\la$, and $\la$ and $\mu$ have the same $h$-content.
\end{enumerate}
This means in particular that canonical basis vectors are weight vectors, so we can consider the canonical basis for a given block.

\begin{eg}
Take $h=5$. Using the rules above, we can compute
\[
f_1f_2f_1f_0^{(3)}f_1f_2f_1f_0\varnothing=(6,4)+q^2(7,3)+q^2(8,2)+q^4(9,1),
\]
so this is the \cbv $\can{6,4}$. We can also compute
\[
f_1f_0f_2f_1f_0^{(2)}f_1f_2f_1f_0\varnothing=(5,3,2)+q^2(5,4,1)+(q^2+1)(6,4)+(q^4+q^2)(7,3)+q^2(8,2)+q^4(9,1).
\]
Subtracting $\can{6,4}$, we obtain
\[
(5,3,2)+q^2(5,4,1)+q^2(6,4)+q^4(7,3),
\]
and this is the \cbv $\can{5,3,2}$.
\end{eg}

One approach to understanding \cbcs is to begin with blocks of small \baw. The case of \baw $0$ is straightforward: given an \hbc $\tau$, the block $\blkw\tau0$ is $1$-dimensional, spanned by the vector $\can\tau=\tau$. In the remainder of this paper we will compute the \cbcs for blocks of \baw $1$ and $2$.

\subsection{Parity}

In the Fock space of type $A^{(1)}$, a result due to Tan \cite{kmt} shows that each canonical basis coefficient is either an even or an odd polynomial; that is either a polynomial in $q^2$, or $q$ times a polynomial in $q^2$. In fact this statement follows from the fact that the canonical basis coefficients coincide with parabolic Kazhdan--Lusztig polynomials, but Tan gives a precise combinatorial criterion for which canonical basis coefficients are even and which are odd. In this section we will prove a corresponding result for type $A^{(2)}$, which appears to be new.

Suppose $\la$ is an \hsp, and recall that $\la'$ denotes the partition conjugate to $\la$. Define the \emph{parity} $\pry\la$ to be the parity of the integer $\sum_{c\gs1}\la'_{ch}$. Our main result is the following.

\begin{propn}\label{party}
Suppose $\la,\mu$ are \hsps with $\mu$ restricted. Then $\dn\la\mu$ is an even function of~$q$ if $\pry\la=\pry\mu$, or an odd function of $q$ if $\pry\la\neq\pry\mu$.
\end{propn}

We will prove this result directly by considering the action of the generators $f_i$ on the Fock space; in contrast, Tan's proof of his result uses the relationship with Kazhdan--Lusztig polynomials, though it seems likely that a direct proof would not be difficult to obtain.

To prove \cref{party}, we define an element $a=\sum_\la a_\la\la\in\scrf$ to be an \emph{even vector} if $a_\la$ is even for all even~$\la$ and odd for all odd $\la$. Alternatively, say that $a$ is an \emph{odd vector} if $a_\la$ is odd for all even~$\la$ and even for all odd $\la$. Say that $a$ is a \emph{parity vector} if it is either an even vector or an odd vector. \cref{party} simply says that each \cbv is a parity vector.

We begin with a lemma.

\begin{lemma}\label{partyf}
Suppose $0\ls i\ls n$ and $\la,\mu$ are \hsps.
\begin{enumerate}
\item
If $i\gs1$, then the coefficient of $\mu$ in $f_i\la$ is an even function of $q$.
\item
The coefficient of $\mu$ in $f_0\la$ is an even function of $q$ if $\pry\la=\pry\mu$, and an odd function of $q$ otherwise.
\item
If $a\in\scrf$ is a parity vector and $k\gs1$, then $f_i^{(k)}a$ is a parity vector.
\end{enumerate}
\end{lemma}

\begin{pf}
(1) is immediate from the definition of the action of $f_i$. To prove (2), we assume $\mu$ appears with non-zero coefficient in $f_0\la$; then $\mu$ is obtained from $\la$ by adding a $0$-node $(r,c)$. Then $\pry\la=\pry\mu$ if $c\equiv1\ppmod h$, while $\pry\la\neq\pry\mu$ if $c\equiv0\ppmod h$.

Examining the coefficient of $\mu$ in $f_0\la$, we can neglect the factor $\prod_{m\in M}f_m$, since this is an even function of $q$. So we just need to show that the number $s$ of addable $0$-nodes of $\la$ to the left of column $c$ minus the number of removable $0$-nodes of $\mu$ to the left of column $c$ is even \iff $c\equiv1\ppmod h$. By considering possible configurations of addable and removable nodes, we find that if $c>1$ then column $1$ contributes $\pm1$ to $s$, while each pair of columns $dh,dh+1$ with $1\ls d<\inp{c/h}$ contributes $0$ or $\pm2$. Finally if $c\equiv1\ppmod h$ then column $c-1$ contributes $\pm1$ to $s$, and we are done.

As a consequence of (1) and (2), we see that if $a\in\scrf$ is a parity vector, then $f_ia$ is a parity vector (of the same parity as $a$). Hence $f_i^ka$ is a parity vector, and hence $f_i^{(k)}a$ is a parity vector, because $f_i^{(k)}a$ is just $f_i^ka$ divided by a function of $q$ which is either even or odd.
\end{pf}

\begin{pf}[Proof of \cref{party}]
The LT algorithm allows us to write
\[
G(\mu)=A(\mu)+\sum_{\substack{\nu\text{ restricted}\\\nu\doms\mu}}b_\nu G(\nu),
\]
where $A(\mu)$ has the form $f_{i_1}^{(k_1)}\dots f_{i_r}^{(k_r)}\varnothing$ and each coefficient $b_\nu$ is symmetric in $q$ and $q^{-1}$. By \cref{partyf} $A(\mu)$ is a parity vector, and by induction we can assume each $G(\nu)$ for $\nu\doms\mu$ is a parity vector. We claim that each $b_\nu$ is a Laurent polynomial in $q$ which is even if $\pry\mu=\pry\nu$, and odd otherwise.

Writing $A(\mu)=\sum a_\la\la$ and examining the coefficient of $\nu$ in the above equation, we obtain
\[
\dn\nu\mu=a_\nu+\sum_{\substack{\xi\text{ restricted}\\\nu\doms\xi\doms\mu}}b_\xi\dn\nu\xi+b_\nu.
\]
The definition of the action of $f_i^{(k)}$ means that $a_\nu$ is a Laurent polynomial; by induction each $b_\xi$ is a Laurent polynomial, and each $\dn\nu\xi$ is a polynomial, so $b_\nu$ is a Laurent polynomial. Now the fact that~$\dn\nu\mu$ is a polynomial divisible by $q$ determines $b_\nu$: if we write
\[
a_\nu+\sum_{\substack{\xi\text{ restricted}\\\nu\doms\xi\doms\mu}}b_\xi\dn\nu\xi=\sum_{t\in\bbz}c_tq^t
\]
with each $c_t\in\bbc$, then
\[
b_\nu=-c_0-\sum_{t<0}c_t(q^t+q^{-t}).
\]
The fact that $A(\mu)$ is a parity vector together with the inductive hypothesis means that $a_\nu+\sum_\xi b_\xi\dn\nu\xi$ is even if $\pry\mu=\pry\nu$ and odd otherwise. Hence the same is true of $b_\nu$.

So by induction our claim about the coefficients $b_\nu$ is proved, and this is enough to show that~$G(\mu)$ is a parity vector.
\end{pf}

\section{Scopes--Kessar pairs}\label{scopeskessarsec}

\subsection{$[w:k]$-pairs}

Our main tool for computing canonical bases for blocks of a given \baw $w$ will be \emph{$[w:k]$}-pairs. These are pairs of blocks $\blkw\sigma w$, $\blkw\tau w$ such that $f_i^{(k)}$ gives a vector space isomorphism $\blkw\sigma w\to\blkw\tau w$ for some $i$; this means that we can deduce the \cbcs for $\blkw\tau w$ from those for $\blkw\sigma w$, and these \cbcs will be very similar (in some cases identical). The genesis of this theory is the fundamental work by Scopes \cite{sc} for blocks of symmetric groups; a version for double covers of symmetric groups was developed by Kessar~\cite{ke} and Kessar--Schaps \cite{kesc}.

First we define a family of involutions $\psi_0,\dots,\psi_n$ on $\hstr$. Take $\la\in\hstr$ and $i\in\{0,\dots,n\}$. Define the \emph{$i$-signature} of $\la$ by working along the edge of $\la$ from left to right, writing a $+$ for each addable $i$-node and a $-$ for each removable $i$-node. Now construct the \emph{reduced $i$-signature} by repeatedly deleting adjacent pairs $+-$ in the $i$-signature. The removable $i$-nodes corresponding to the $-$ signs in the reduced $i$-signature are called \emph{normal} $i$-nodes of $\la$, while the addable $i$-nodes corresponding to the $+$ signs are called the \emph{conormal} nodes of $\la$.

Suppose $\la$ has $r$ normal $i$-nodes and $s$ conormal $i$-nodes. Define $\psi_i(\la)$ by:
\begin{itemize}
\item
adding the leftmost $s-r$ conormal $i$-nodes, if $s\gs r$;
\item
removing the rightmost $r-s$ normal $i$-nodes, if $r\gs s$.
\end{itemize}
It is an easy exercise to check that $\psi_i(\la)$ is also an \hsp.

\begin{rmk}
The involution $\psi_i$ derives from the crystal structure of $\scrf$: the theory of crystal bases (see \cite{hoka} for an introduction) defines a directed graph (the \emph{crystal} of $\scrf$) with vertex set $\hstr$ and edges labelled by residues $i\in\{0,\dots,n\}$; the subgraph formed by the edges labelled $i$ is just a disjoint union of directed paths, and the effect of $\psi_i$ is to reverse each of these paths.
\end{rmk}

\begin{eg}
Suppose $h=3$, $\la=(5,4,2,1)$ and $i=0$. The addable and removable $0$-nodes of $\la$ are indicated in the following diagram.
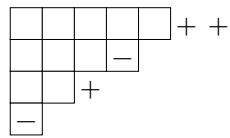
\[
\gyoung(;;;;;:+:+,;;;;-,;;:+,;-)
\]
So the $0$-signature of $\la$ is $-+-++$, so that the reduced $0$-signature is $-++$; $\la$ has one normal $0$-node $(4,1)$, and two conormal $0$-nodes $(1,6)$ and $(1,7)$. So $\psi_0(\la)=(6,4,2,1)$.
\end{eg}

We remark on two very easy special cases.

\begin{lemma}\label{psispec}
Suppose $\la\in\hstr$ and $i\in\{0,\dots,n\}$.
\begin{enumerate}
\item
If $\la$ has no removable $i$-nodes, then $\psi_i(\la)$ is obtained by adding all the addable $i$-nodes to $\la$.
\item
If $\la$ has no addable $i$-nodes, then $\psi_i(\la)$ is obtained by removing all the removable $i$-nodes from~$\la$.
\end{enumerate}
\end{lemma}

We also note some properties of the functions $\psi_i$, which are well known and easy to prove. (The first of these comes from the following simple observation: if $\la,\mu\in\hstr$ and $\mu$ is obtained by adding an $i$-node $(r,c)$ to $\la$, then $(r,c)$ is the rightmost normal $i$-node of $\mu$ \iff it is the leftmost conormal $i$-node of $\la$.)

\begin{lemma}\label{psibasic}
Suppose $\la\in\hstr$ and $i\in\{0,\dots,n\}$.
\begin{enumerate}
\item
$\psi_i^2(\la)=\la$.
\item
$\psi_i(\la)$ is restricted \iff $\la$ is restricted.
\end{enumerate}
\end{lemma}

Now we consider the particular case of \hbcs.

\begin{propn}\label{corerefl}
Suppose $\sigma$ is an \hbc and $0\ls i\ls n$. Then:
\begin{enumerate}[ref=(\arabic*)]
\item\label{noaddrem}
$\sigma$ cannot have both addable and removable $i$-nodes;
\item\label{psihb}
$\psi_i(\sigma)$ is an \hbc;
\item\label{hbpsi}
if $\tau$ is an \hbc obtained by adding some $i$-nodes to $\sigma$, then $\tau=\psi_i(\sigma)$.
\end{enumerate}
\end{propn}

\begin{pf}
We use the abacus, in particular \cref{addnodesabacus,coreflush}.

First suppose $1\ls i<n$. If $\sigma$ has addable $i$-nodes, then there is at least one bead on runner $i$ of the \abd for $\sigma$ with an empty space immediately to its right. But by \cref{coreflush} the assumption that $\sigma$ is an \hbc means that every bead in the \abd has a bead immediately above it. So runners $i$ and $i+1$ of the \abd have the following form.
\[
\abacus(bb,bn,bn,vv,bn,bn,nn)
\]
So $\sigma$ has no removable $i$-nodes, proving \ref{noaddrem}. Moreover, the \abd for $\psi_i(\sigma)$ is obtained by simply switching runners $i$ and $i+1$ (and also runners $-i-1$ and $-i$) so every bead in the \abd for $\psi_i(\sigma)$ has a bead immediately above it. So by \cref{coreflush} $\psi_i(\sigma)$ is an \hbc, so \ref{psihb} holds. For \ref{hbpsi}, apply \ref{noaddrem} to $\tau$: since by assumption $\tau$ has removable $i$-nodes, it cannot have addable $i$-nodes, so in order to obtain $\tau$ from $\sigma$, all the addable $i$-nodes must have been added, and therefore $\tau=\psi_i(\sigma)$.

The case $i=n$ is very similar, except that here there is only one pair of runners to consider, namely~$n$ and $-n$, and the phrase ``to the right of'' must be reinterpreted appropriately.

The case $i=0$ is also similar, except that here there are three runners to consider, namely runners $-1$, $0$ and $1$. Here the fact that every bead in the \abd for $\sigma$ has a bead immediately above it and the assumption that $\sigma$ has at least one addable $0$-node, together with the symmetry of the abacus, means that the configuration on runners $-1$, $0$ and $1$ of the \abd for $\sigma$ is as follows.
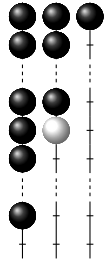
\[
\abacus(bbb,bbn,vvv,bbn,bxn,bnn,vvv,bnn,nnn)
\]
Now the \abd for $\psi_0(\sigma)$ is obtained by switching runners $-1$ and $1$, and the proof works as for $i\gs1$.
\end{pf}

A by-product of the above proof is that if an \hbc $\sigma$ has addable $0$-nodes, it must have an odd number of them; we will use this observation repeatedly.

Now we return to arbitrary \hsps.

\begin{lemma}\label{phicorewt}
Suppose $0\ls i\ls n$ and $\la\in\hstr$, and that $\la$ has \hbc $\sigma$ and \hbw $w$. Then~$\psi_i(\la)$ has \hbc $\psi_i(\sigma)$ and \hbw $w$.
\end{lemma}

\begin{pf}
We assume $1\ls i<n$, with the cases $i=0$ and $i=n$ being similar. Let $\delta_i(\la)$ equal the number of addable $i$-nodes of $\la$ minus the number of removable $i$-nodes of $\la$. First we claim that $\delta_i(\la)=\delta_i(\sigma)$. To see this, we use the \abd for $\la$. \cref{addnodesabacus} shows that, if we take $N\gg0$, then $\delta_i(\la)$ equals the number of beads on runner $i$ of the \abd after position $i-Nh$ minus the number of beads on runner $i+1$ after position $i+1-Nh$. Constructing the \abd for $\sigma$ involves moving beads up their runners, so does not affect these numbers of beads; so $\delta_i(\sigma)=\delta_i(\la)$.

The way the reduced $i$-signature is constructed means that $\delta_i(\la)$ is also the number of conormal $i$-nodes of $\la$ minus the number of normal $i$-nodes. So $\psi_i(\la)$ is obtained from $\la$ by adding $\delta_i(\la)$ $i$-nodes if $\delta_i(\la)\gs0$, or removing $-\delta_i(\la)$ $i$-nodes if $\delta_i(\la)<0$. Since $\delta_i(\la)=\delta_i(\sigma)$, the same applies for~$\sigma$ and~$\psi_i(\sigma)$.

We can easily compare the $h$-contents of $\la$ and $\sigma$: removing an $h$-bar entails removing two $k$-nodes for each $k\in\{0,\dots,n-1\}$ and one $n$-node; so the $h$-content of $\sigma$ is obtained from the $h$-content of $\la$ by removing $w$ copies of $n$ and $2w$ copies of $k$ for each $0\ls k<n$. The previous paragraph implies that the same relationship holds between the $h$-contents of $\psi_i(\la)$ and $\psi_i(\sigma)$. So $\psi_i(\la)$ has the same $h$-content as an \hsp with \hbc $\psi_i(\sigma)$ and \hbw $w$, and so by \cref{contcore}~$\psi_i(\la)$ has \hbc $\psi_i(\sigma)$ and \hbw $w$.
\end{pf}

Now we can introduce Scopes--Kessar pairs. Suppose $\sigma$ is an \hbc with $k$ addable $i$-nodes, where $k\gs1$. Let $\tau=\psi_i(\sigma)$. We say that $\blkw\sigma w$ and $\blkw\tau w$ form a \emph{$[w:k]$-pair} of residue~$i$. It follows from \cref{phicorewt} that $\psi_i$ restricts to a bijection between $\pbw\sigma w$ and $\pbw\tau w$.

We define a partition $\la\in\pbw\sigma w$ to be \emph{unexceptional} (for the pair $(\blkw\sigma w,\blkw\tau w)$) if it has no removable $i$-nodes, and we define a partition $\mu\in\pbw\tau w$ to be unexceptional if it has no addable $i$-nodes. We consider the effect of the operator $f_i^{(k)}$.

\begin{propn}\label{unex}
Suppose $\sigma,\tau,k,i$ are as above.
\begin{enumerate}[ref=\cref{unex}(\arabic*)]
\item\label{fun}
If $\la\in\pbw\sigma w$ is unexceptional, then $\la$ has exactly $k$ addable $i$-nodes, and $f_i^{(k)}\la=\psi_i(\la)$.
\item\label{fgun}
If $\mu\in\pbw\sigma w$ is restricted and $\can\mu$ is a linear combination of unexceptional partitions, then (extending $\psi_i$ linearly)
\[
\can{\psi_i(\mu)}=\psi_i(\can\mu).
\]
\end{enumerate}
\end{propn}

\begin{pfenum}
\item
For the first statement we use the fact (shown in the proof of \cref{phicorewt}) that the number of addable $i$-nodes minus the number of removable $i$-nodes is the same for $\la$ as it is for $\sigma$; by \cref{corerefl} $\sigma$ has $k$ addable $i$-nodes and no removable $i$-nodes, so the same is true for $\la$.

Now the second statement follows from the formula for the action of $f_i^{(k)}$.
\item
Write $\can\mu=\sum_\la \dn\la\mu\la$, with each $\la$ unexceptional. Then by (1)
\[
f_i^{(k)}\can\mu=\sum_\la \dn\la\mu\psi_i(\la).
\]
Since this vector is bar-invariant and each coefficient is divisible by $q$ except for the coefficient of $\psi_i(\mu)$ which equals $1$, this vector must equal $\can{\psi_i(\mu)}$.\qedhere
\end{pfenum}

\cref{unex} provides a result analogous to the \emph{Scopes--Kessar equivalences} for blocks of symmetric groups and their double covers; it says that if every partition in $\blkw\sigma w$ is unexceptional, then~$\blkw\sigma w$ and $\blkw\tau w$ have the same canonical basis (up to relabelling of basis elements); we say that $\blkw\sigma w$ and~$\blkw\tau w$ are \emph{Scopes--Kessar equivalent} in this case. (This is essentially the same as saying that $(\sigma,\tau)$ is a \emph{$w$-compatible pair}, as defined by Kessar and Schaps \cite[Definition~3.1]{kesc}.)

In fact it is possible to say exactly (in terms of $w$, $k$ and $i$) when we have a Scopes--Kessar equivalence. The next result is essentially \cite[Lemma 3.8]{lesc}, but we provide a self-contained statement and proof to avoid translating notation from one setting to another, and because our abacus convention affords a shorter proof than the one in \cite{lesc}.

\begin{thm}\label{skequiv}
Suppose $\blkw\si w$ and $\blkw\tau w$ form a $[w:k]$-pair of residue $i$. Then $\blkw\si w$ and $\blkw\tau w$ are Scopes--Kessar equivalent \iff
\[
w\ls
\begin{cases}
\frac12(k-1)&\text{if }i=0
\\
k&\text{if }1\ls i<n
\\
2k+1&\text{if }i=n.
\end{cases}
\]
\end{thm}

\begin{pf}
We consider the three cases separately.
\begin{description}
\item
$1\ls i<n$\\
This case is the most closely related to Scopes's original equivalence for the symmetric group. Suppose that in the \abd for $\si$ the lowest bead on runner $i$ is in position $ch+i$. Then the lowest bead on runner $i+1$ is in position $(c-k)h+i+1$. If $\la\in\blkw\si w$, then the \abd of $\la$ is obtained from the \abd of $\si$ by moving beads down their runners. Suppose that $a$ of these moves take place on runner $i$ and $b$ on runner $i+1$. Then $a+b\ls w$. The lowest bead on runner $i+1$ of the \abd for $\la$ is in position $(c-k+b)h+i+1$ or higher, since each bead move can only move the lowest bead down by at most one row. Similarly, the highest empty position on runner $i$ is position $(c+1-a)h+i$ or lower. If $\la$ has a removable $i$-node, then there is a bead on runner $i+1$ with an empty position immediately to its left, which means that $c-k+b\gs c+1-a$, and hence that $w\gs a+b\gs k+1$. So if $w\ls k$ then there are no exceptional partitions in $\blkw\si w$. Conversely, if $w>k$, then we can easily construct an exceptional partition: starting from the \abd for $\si$, we move the lowest bead on runner $i$ down $w$ positions, and the lowest $w$ beads on runner $-i$ down one position each.
\item
$i=0$\\
In this case the lowest bead on runner $-1$ of the \abd for $\si$ is in position $\frac12(k-1)-1$. Take $\la\in\blkw\si w$, and suppose that in constructing the \abd of $\la$ from $\si$ we make $a$ bead moves on runner $-1$ and $b$ moves on runner $0$. As in the previous case, we find that it is possible for $\la$ to have a removable $0$-node (i.e.\ for there to be a bead on runner $0$ with an empty space immediately to its left) \iff $a+b>\frac12(k-1)$.
\item
$i=n$\\
In this case the lowest bead on runner $n$ is in position $(k-1)h+n$. Suppose that when we construct the \abd of $\la\in\blkw\si w$ from $\si$, we make $a$ bead moves on runner $n$. If $\la$ has a removable $n$-node, then there is a bead at position $ch+n+1$ for some $c$ such that position $ch+n$ is unoccupied. From the symmetry of the abacus, position $(-c-1)h+n$ is also unoccupied. By swapping $c$ and $-c-1$ if necessary, we can assume $c\gs0$. Now to construct the \abd for $\la$ from the \abd for $\si$, the lowest $k+c+1$ beads on runner $n$ must be moved down (to create the space at position $(-c-1)h+n$); in addition, if $c\ls k$, then the lowest $k-c+1$ beads must move down again to create the additional space at position $ch+n$. So
\[
w\gs a\gs
\begin{cases}
k+c+1&\text{if }c>k
\\
(k+c+1)+(k-c+1)&\text{if }c\ls k,
\end{cases}
\]
which in either case is greater than $2k+1$.

Conversely, if $w\gs2k+2$ then we can construct an exceptional partition: move the lowest $2k+2$ beads on runner $n$ down one row, and the lowest bead on runner $-n$ down $2k+2$ rows.\qedhere
\end{description}
\end{pf}

\subsection{$[w:k]$-pairs and orders on partitions}

We end this section by considering how the various orderings defined in \cref{hspsec} change as we pass through a $[w:k]$-pair. Given $\la,\mu\in\hstr$ and $i\in\{0,\dots,n\}$, we write $\link\la{\hat\la}i$ if $\hat\la$ can be obtained from $\la$ by adding and/or removing $i$-nodes. Observe that if $\blkw\sigma w$ and $\blkw\tau w$ form a $[w:k]$-pair of residue $i$ and $\la\in\pbw\sigma w$, $\hat\la\in\pbw\tau w$ with either $\la$ or $\hat\la$ unexceptional, then $\link\la{\hat\la}i$ \iff $\hat\la=\psi_i(\la)$.

\begin{lemma}\label{unexlex}
Suppose $\blkw\sigma w$ and $\blkw\tau w$ form a $[w:k]$-pair of residue $i$. Suppose $\la,\mu\in\pbw\sigma w$ and $\hat\la,\hat\mu\in\pbw\tau w$ with $\link\la{\hat\la}i$ and $\link\mu{\hat\mu}i$, and that $\la$ is unexceptional.
\begin{enumerate}[ref=(\arabic*)]
\item\label{unexlexlex}
If $\la\glex\mu$, then $\hat\la\glex\hat\mu$.
\item
If $\la\lcolex\mu$, then $\hat\la\lcolex\hat\mu$.
\end{enumerate}
\end{lemma}

\begin{pf}
We prove the first part only, as the proof of the second part is very similar.

Since $\la\glex\mu$, there is $r\gs1$ such that $\la_r>\mu_r$ while $\la_s=\mu_s$ for all $s<r$. So for any $s<r$, the number of addable nodes in row $s$ is the same for $\la$ as for $\mu$. Since $\la$ is unexceptional, $\hat\la$ is obtained from $\la$ by adding all the addable $i$-nodes, so we get $\hat\la_s\gs\hat\mu_s$ for all $s<r$. If $\hat\la_s>\hat\mu_s$ for any $s\ls r$ then we are done, so assume $\hat\la_s=\hat\mu_s$ for $s<r$, and $\hat\la_r\ls\hat\mu_r$. This means that $\mu$ has at least one addable $i$-node in row $r$, and that the node $(r,\la_r)$ is an $i$-node. This node cannot be a removable $i$-node of $\la$ because $\la$ is unexceptional, which means that $\la_{r+1}=\la_r-1$ and $\la_r\nequiv1\ppmod h$. Furthermore, we claim that $\hat\la_r=\hat\mu_r$: the only way this could in theory fail is if $i=0$ and $\hat\mu_r=\mu_r+2=\la_r+1$, but then the node $(r,\mu_r+2)$ would also be an addable node of $\la$, so that $\hat\la_r=\hat\mu_r$.
 
The deduction that $\la_r\nequiv1\ppmod h$ means that $\mu_r\nequiv0\ppmod h$, and hence $\mu_{r+1}<\mu_r$. Hence $\hat\mu_{r+1}<\mu_r\ls\la_r-1=\la_{r+1}\ls\hat\la_{r+1}$, and so $\hat\la\glex\hat\mu$.
\end{pf}

\section{Blocks of \baw $1$}\label{wt1sec}

In this section we determine the $q$-decomposition numbers for blocks of \baw $1$. The decomposition numbers for spin blocks of symmetric groups of defect $1$ are known, thanks to M\"uller \cite{muller}; the results in the present section provide a $q$-analogue of M\"uller's results.

Throughout this section we fix an \hbc $\tau$. Recall that we define $\bb\tau$ to be the number of parts of $\tau$ less than $h$.

\begin{propn}\label{ptnsinwt1}
There are exactly $n+1$ partitions in $\pbw\tau1$. These are totally ordered by the dominance order, with all except the most dominant one being restricted.
\end{propn}

\begin{pf}
There are three different ways to add an $h$-bar to $\tau$ to obtain an \hsp, for which we use M\"uller's notation.
\begin{description}
\item[type $\oplus$:]
Given $a\in\tau$ such that $a>0$ and $a+h\notin\tau$, define $\tau[a+h]=\tau\sqcup(a+h)\setminus(a)$. This partition is restricted unless $a=\tau_1$.
\item[type $0$:]
Define $\tau[h]=\tau\sqcup(h)$. This partition is restricted unless $\tau=\varnothing$.
\item[type $\ominus$:]
Given $b\in\{n+1,\dots,h-1\}$ such that neither $b$ nor $h-b$ lies in $\tau$, define $\tau[b]=\tau\sqcup(b,h-b)$. This partition is restricted.
\end{description}
Now observe that if $\tau[i]$ and $\tau[j]$ are two partitions constructed in this way, then $\tau[i]\domsby\tau[j]$ \iff $i<j$. Moreover, the number of partitions of type $\oplus$ equals $\bb\tau$, and (because $\tau$ cannot contain integers $a$ and $h-a$) the number of partitions of type $\ominus$ equals $n-\bb\tau$. So there are $n+1$ partitions in $\pbw\tau1$ altogether, and only the most dominant fails to be restricted.
\end{pf}

In view of \cref{ptnsinwt1}, we can label the partitions in $\pbw\tau1$ as $\tau(0)\domsby\cdots\domsby\tau(n)$. Now we can give the main result of this section.

\begin{thm}\label{mainwt1}
Suppose $\tau$ is an \hbc, and let $\tau(0)\domsby\cdots\domsby\tau(n)$ be the partitions in $\pbw\tau1$. Then
\[
\dn{\tau(r)}{\tau(s)}=
\begin{cases}
1&\text{if }r=s\\
q&\text{if }r=s+1\text{ and }h\in\tau(r)\\
q^2&\text{if }r=s+1\text{ and }h\notin\tau(r)\\
0&\text{otherwise}.
\end{cases}
\]
\end{thm}

For example, with $h=7$ and $\tau=(4,2)$, the $q$-decomposition numbers are given by the following matrix. (In all explicit matrices in this paper, we use $\cdot$ to mean $0$.)
\[
\begin{array}{c|ccc}
&\rt{6,4,2,1}&\rt{7,4,2}&\rt{9,4}\\\hline
(6,4,2,1)&1&\cdot&\cdot\\
(7,4,2)&q&1&\cdot\\
(9,4)&\cdot&q^2&1\\
(11,2)&\cdot&\cdot&q^2
\end{array}
\]

The exceptional entry $q$ in the $q$-decomposition matrix reflects the change in parity as we we read the partitions in $\pbw\tau1$: the partitions of type $\oplus$ have the same parity as $\tau$, while those of type $0$ or $\ominus$ have the opposite parity.

We prove \cref{mainwt1} by induction. For the initial cases of the induction, we assume that $\tau$ has the form $(l,l-1,\dots,1)$ for some $0\ls l\ls n$. In this situation, \cref{mainwt1} can be re-cast as follows.

\begin{propn}\label{mainwt12core}
Suppose $\tau=(l,l-1,\dots,1)$ for $0\ls l\ls n$, and suppose $\mu\in\pbw\tau1$ is restricted. Then one of the following occurs.
\begin{enumerate}
\item
$\mu=\tau\sqcup(h-b,b)$ for $l+2\ls b\ls n$, and $\can\mu=\mu+q^2\bigl(\tau\sqcup(h-b+1,b-1)\bigr)$.
\item
$\mu=\tau\sqcup(h-l-1,l+1)$ and $l\ls n-1$, and $\can\mu=\mu+q\bigr(\tau\sqcup(h)\bigr)$.
\item
$\mu=\tau\sqcup(h)$ and $l\gs1$, and $\can\mu=\mu+q^2\bigl(\tau\sqcup(h+1)\setminus(1)\bigr)$.
\item
$\mu=\tau\sqcup(a+h)\setminus(a)$ with $1\ls a\ls l-1$, and $\can\mu=\mu+q^2\bigl(\tau\sqcup(a+h+1)\setminus(a+1)\bigr)$.
\end{enumerate}
\end{propn}

\begin{pf}
It is clear that $\mu$ satisfies one (and only one) of the given conditions, so we just need to calculate~$\can\mu$. In each case we construct an \hbc $\xi$, and show that the given vector can be constructed from $\can\xi=\xi$ by applying the operators $f_i$. The defining properties of the canonical basis then guarantee that this vector must equal~$\can\mu$.
\begin{enumerate}
\item
Let $\xi=\tau\sqcup(h-b,b-1)$. Then
\begin{align*}
f_{b-1}\can\xi&=\mu+q^2\bigl(\tau\sqcup(h-b+1,b-1)\bigr).
\\
\intertext{\item
Let $\xi=\tau\sqcup(h-l-1)$. Then}
f_0f_1\dots f_l\can\xi&=\mu+q\bigr(\tau\sqcup(h)\bigr).
\\
\intertext{\item
Let $\xi=\tau\sqcup(h-1)\setminus(1)$. Then}
f_0^{(2)}\can\xi&=\mu+q^2\bigl(\tau\sqcup(h+1)\setminus(1)\bigr).
\\
\intertext{\item
Let $\xi=\tau\sqcup(a+h)\setminus(a+1)$. Then}
f_a\can\xi&=\mu+q^2\bigl(\tau\sqcup(a+h+1)\setminus(a+1)\bigr).\qedhere
\end{align*}
\end{enumerate}
\end{pf}

\begin{pf}[Proof of \cref{mainwt1}]
We proceed by induction on $\card\tau$. Let $l=\len\tau$, and consider the three possibilities in \cref{non2core}. If $\tau=(l,l-1,\dots,1)$ with $l\ls n$, then \cref{mainwt12core} gives the result. Alternatively, there is a residue $i$ such that either $i\neq0$ and $\tau$ has a removable $i$-node, or $i=0$ and $\tau$ has at least three removable $i$-nodes. So define the \hbc $\xi$ by removing all the removable $i$-nodes from~$\tau$. Then $\blkw\xi1$ and $\blkw\tau1$ are Scopes--Kessar equivalent, by \cref{skequiv}. Since the lexicographic order refines the dominance order, we have
\[
\xi(0)\llex\cdots\llex\xi(n),\qquad\tau(0)\llex\cdots\llex\tau(n),
\]
and so by \cref{unexlex} $\psi_i(\xi(r))=\tau(r)$ for all $r$. So by \cref{unex} $\dn{\xi(r)}{\xi(s)}=\dn{\tau(r)}{\tau(s)}$ for all $r,s$. Finally, the Scopes--Kessar equivalence sends the partition $\xi\sqcup(h)$ to $\tau\sqcup(h)$, and the result follows from the inductive hypothesis.
\end{pf}

\section{Blocks of \baw $2$}\label{wt2mainsec}

Now we come to the main object of study in this paper: blocks of \baw $2$. In this section we further develop the combinatorics of blocks of \baw $2$ and state our main theorem. We continue to work with a fixed \hbc $\tau$.

\subsection{Abacus notation}

Our work with blocks of \baw $2$ requires some uniform notation. Suppose $\la\in\pb\tau$. First we define two integers which we call the \emph{bar positions} of $\la$. Each time we remove an $h$-bar, we move beads on the abacus from positions $a,h-a$ to positions $a-h,-a$ for some $a\gs n+1$. We call $a$ the bar position of this $h$-bar. Define $\lo\la$ and $\hi\la$ to be the two bar positions of $\la$, with $\lo\la\ls\hi\la$.

We also define a notation for $\la$ based on the ``active'' runners of its \abd, i.e.\ those runners in which not every bead has a bead immediately above it. This is analogous to the notation introduced by Scopes \cite{sc2} in weight $2$ blocks of symmetric groups, and is defined as follows.
\begin{itemize}
\item
If (the \abd for) $\la$ is obtained from $\tau$ by moving the lowest bead down one space each on runners $-j<-i<i<j$, then we write $\la=\ba{i,j}$ (or $\ba{j,i}$).
\item
If $\la$ is obtained from $\tau$ by moving the lowest bead down one space on runners $-j<j$ and moving a bead from $-h$ to $h$, then we write $\la=\ba{0,j}$.
\item
If $\la$ is obtained from $\tau$ by moving the lowest bead on runner $i$ and the second-lowest bead on runner $-i$ down two spaces each (where $i\neq0$), then we write $\la=\ba i$.
\item
The partition $\la=\tau\sqcup(2h)$ is written as $\ba0$.
\item
The partition $\la=\tau\sqcup(h^2)$ is written as $\ba{0,0}$.
\end{itemize}

\begin{eg}
Take $h=7$. The partitions $\la=(15,9,2)$ and $\mu=(9,8,6,2,1)$ both have $7$-bar-core $(9,2,1)$, and satisfy $(\lo\la,\hi\la)=(8,15)$ and $(\lo\mu,\hi\mu)=(6,8)$. The abacus notations for $\la$ and $\mu$ are $\ba1$ and $\ba{-1}$ respectively.
\[
\begin{array}{c@{\qquad}c}
\abacus(bbnbbbb,bnbbbbb,bnbxnbn,nnnnnbn,nnnnbnn)
&
\abacus(bbbbbbb,bnnbnbb,bnnxbbn,nnbnbbn,nnnnnnn)
\\[24pt]
\la=(15,9,2)
&
\mu=(9,8,6,2,1)
\end{array}
\]
\end{eg}

\subsection{The dominance order}\label{domdef2sec}

Now we consider the dominance order in $\pb\tau$, which will be central to the our formula for the \cbcs. Given $\la,\mu\in\pb\tau$, we give a simple sufficient criterion for $\la\domby\mu$ in terms of the bar positions $\lo\la,\hi\la$ introduced above. This is an analogue of a lemma in Richards's paper \cite[Lemma~4.4]{rich}.

\begin{propn}\label{domlohi}
Suppose $\la,\mu\in\pb\tau$. If $\lo\la\ls\lo\mu$ and $\hi\la\ls\hi\mu$, then $\la\domby\mu$.
\end{propn}

\begin{pf}
Replacing an integer $a$ with $a+h$ in a strict partition entails adding one node to each of columns $a+1,\dots,a+h$; similarly, inserting the parts $a,h-a$ with $a\ls n$ entails adding two nodes to each of columns $1,\dots,a$, and one node to each of columns $a+1,\dots,h-a$.

So if $\lo\la\ls\lo\mu$ and $\hi\la\ls\hi\mu$, then when we construct $\la$ from $\tau$, we add nodes in earlier columns than when we construct $\mu$ from $\tau$; so $\la'\dom\mu'$, and therefore $\la\domby\mu$ by \cref{conjdom}.
\end{pf}

\begin{rmk}
In fact Richards's result is stronger than ours in that the condition he gives is necessary and sufficient for $\la\domby\mu$. This is not quite true for \cref{domlohi}; for example, take $h=7$, $\la=(6,5,2,1)$ and $\mu=(7,4,3)$. Then $\la\domby\mu$ even though $\lo\la=5$ and $\lo\mu=4$. A necessary and sufficient condition for $\la\domby\mu$ in terms of $\lo\la,\hi\la,\lo\mu,\hi\mu$ would not be hard to obtain, but we do not need it here.
\end{rmk}

\subsection{Leg lengths}

Richards's formula for decomposition numbers of defect $2$ blocks of symmetric groups involves leg lengths of rim hooks. Here we introduce the corresponding combinatorics for our situation.

Suppose $\la\in\pb\tau$, and let $\lo\la,\hi\la$ be bar positions for $\la$, as defined above. Recall that we write $\ope xy{\la\cap\tau}$ for the number of integers strictly between $x$ and $y$ that lie in both $\la$ and $\tau$. Take $\lohi\in\{\lo\la,\hi\la\}$, and define the \emph{leg length} corresponding to $\lohi$ as follows:
\begin{itemize}
\item
if $\lohi\gs h$, define the leg length to be $\ope{\lohi-h}\lohi{\la\cap\tau}$;
\item
if $n<\lohi<h$, define the leg length to be $h-\lohi+\ope{h-\lohi}\lohi{\la\cap\tau}$.
\end{itemize}

We define $\ddd\la$ to be the absolute value of the difference between the leg lengths of $\la$.

\begin{rmk}
Leg lengths for $h$-bars are also defined by Hoffman and Humphreys in \cite[p.185]{hohum}. However, our definition differs slightly from theirs. For example, take $h=7$ and $\la=(10,5,4)$. Then $\lo\la=4$ and $\hi\la=10$, giving leg lengths $1$ and $3$, and hence $\ddd\la=2$. Using the notion of leg length from \cite{hohum}, the leg length corresponding to $\lohi=4$ would be $2$, giving $\ddd\la=1$.
\end{rmk}

The key to Richards's combinatorial formula for decomposition numbers is the interplay between the $\ddd$ function and the dominance order. In our setting, the following result will be important.

\begin{propn}\label{domddd}
Suppose $\la,\mu\in\pb\tau$. If $\la$ and $\mu$ are incomparable in the dominance order, then $|\ddd\la-\ddd\mu|\gs2$.
\end{propn}

\begin{pf}
By \cref{domlohi}, we may assume (interchanging $\la$ and $\mu$ if necessary) that $\hi\la>\hi\mu$ and $\lo\la<\lo\mu$. Then we claim that $\ddd\la\gs\ddd\mu+2$.

There are several cases to check, depending on the relative order of the eight integers $\hi\la$, $\hi\mu$, $\lo\mu$, $\lo\la$, $\hi\la-h$, $\hi\mu-h$, $\lo\mu-h$, $\lo\la-h$. We show the calculations in three example cases, leaving the reader to check the other cases.
\begin{enumerate}
\item
First suppose $\hi\la<h$. Then the leg lengths of $\la$ are (in decreasing order)
\[
h-\lo\la+\ope{h-\lo\la}{\lo\la}\tau,\qquad h-\hi\la+\ope{h-\hi\la}{\hi\la}\tau
\]
and similarly for $\mu$, giving
\[
\ddd\la-\ddd\mu=\hi\la-\hi\mu+\lo\mu-\lo\la-\ope{h-\hi\la}{h-\hi\mu}\tau-\ope{h-\lo\mu}{h-\lo\la}\tau-\ope{\lo\la}{\lo\mu}\tau-\ope{\hi\mu}{\hi\la}\tau.
\]
Now the fact that $\tau$ does not contain two parts summing to $h$ implies that
\begin{align*}
\ope{h-\hi\la}{h-\hi\mu}\tau+\ope{\hi\mu}{\hi\la}\tau&<\hi\la-\hi\mu\\
\intertext{and}
\ope{h-\lo\mu}{h-\lo\la}\tau+\ope{\lo\la}{\lo\mu}\tau&<\lo\mu-\lo\la,\\
\end{align*}
so that $\ddd\la-\ddd\mu\gs2$.
\item
Now suppose $\lo\mu=\hi\mu=h$. If $\hi\la-h\gs h-\lo\la$, then we can obtain $\mu$ from $\la$ by moving nodes from columns $\hi\la-h+1,\dots,\hi\la-h+\lo\la$ to columns $h-\lo\la+1,\dots,h$ and from columns $\hi\la-h+\lo\la+1,\dots,\hi\la$ to columns $\lo\la+1,\dots,h$; each moved node moves to an earlier column, which gives $\la\dom\mu$, contrary to assumption.

So instead we must have $\hi\la-h<h-\lo\la<\lo\la<h<\hi\la$. Now the leg lengths for $\la$ are
\[
h-\lo\la+\ope{h-\lo\la}{\lo\la}\tau,\qquad\ope{\hi\la-h}{\hi\la}\tau
\]
while obviously $\ddd\mu=0$, so that
\begin{align*}
\ddd\la-\ddd\mu&=h-\lo\la-\ope{\hi\la-h}{h-\lo\la}\tau-\ope{\lo\la}h\tau-\ope h{\hi\la}\tau.
\\
\intertext{Because $\tau$ is an \hbc, $\ope h{\hi\la}\tau\ls\ope0{\hi\la-h}\tau$, so that}
\ddd\la-\ddd\mu&\gs h-\lo\la-\ope0{\hi\la-h}\tau-\ope{\hi\la-h}{h-\lo\la}\tau-\ope{\lo\la}h\tau\\
&=h-\lo\la-\ope0{h-\lo\la}\tau-\ope{\lo\la}h\tau+1,
\end{align*}
with the $1$ arising because $\hi\la-h\in\tau$. But now the fact that $\tau$ does not contain two parts summing to $h$ gives $\ope0{h-\lo\la}\tau+\ope{\lo\la}h\tau\ls h-\lo\la-1$, so that $\ddd\la-\ddd\mu\gs2$.
\item
Now suppose $\lo\la>h$ and $\hi\la<\lo\la+h$. Then the leg lengths for $\la$ are
\[
\ope{\lo\la-h}{\lo\la}\tau-1,\qquad\ope{\hi\la-h}{\hi\la}\tau,
\]
with the $-1$ occurring because $\hi\la-h$ lies between $\lo\la-h$ and $\lo\la$ and belongs to $\tau$ but not $\la$, so does not get counted in the leg length. A similar statement applies for $\mu$, so that
\begin{align*}
\ddd\la-\ddd\mu&=\ope{\lo\la-h}{\lo\la}\tau-\ope{\lo\mu-h}{\lo\mu}\tau+\ope{\hi\mu-h}{\hi\mu}\tau-\ope{\hi\la-h}{\hi\la}\tau\\
&=\ope{\lo\la-h}{\lo\mu-h}\tau+\ope{\hi\mu-h}{\hi\la}\tau-\ope{\lo\la}{\lo\mu}\tau-\ope{\hi\mu}{\hi\la}\tau+2
\end{align*}
with the $2$ arising from counting the occurrences of $\lo\mu-h,\hi\mu-h,\hi\la-h$, which all belong to $\tau$. But $\tau$ is an \hbc, so
\[
\ope{\lo\la-h}{\lo\mu-h}\tau\gs\ope{\lo\la}{\lo\mu}\tau,\qquad\ope{\hi\mu-h}{\hi\la-h}\tau\gs\ope{\hi\mu}{\hi\la}\tau,
\]
giving $\ddd\la-\ddd\mu\gs2$.\qedhere
\end{enumerate}
\end{pf}

Now (continuing the analogy with Richards's work) we define the \emph{colour} of a partition $\la\in\pb\tau$. The way we define this depends on the value of $\ddd\la$.
\begin{description}
\item
$\ddd\la=0$\\
Suppose first that $\la$ has a $2h$-bar, i.e.\ either $\lo\la=\hi\la-h$ or $\lo\la+\hi\la=2h$. We define the leg length of this $2h$-bar analogously to the leg length of an $h$-bar: if $\hi\la\gs2h$ then we define the leg length as $l=\card{\la\cap\tau\cap\{\hi\la-2h+1,\dots,\hi\la-1\}}$, while if $\hi\la<2h$ we define the leg length to be $l=2h-\hi\la+\card{\la\cap\tau\cap\{2h-\hi\la+1,\dots,\hi\la-1\}}$. Now we say $\la$ is \emph{black} if $l+2\bb\tau$ is congruent to $0$ or $3$ modulo $4$, and \emph{white} otherwise.

Alternatively, suppose $\la$ has two $h$-bars, with common leg length $l$. Say that $\la$ is \emph{black} if $l+\bb\tau$ is odd, and \emph{white} otherwise.
\item
$\ddd\la=1$\\
If $\hi\la>h$, we say that $\la$ is \emph{grey}. If $\hi\la\ls h$, then let $l,l+1$ be the two leg lengths of $\la$. We say that~$\la$ is \emph{black} if $l+\bb\tau$ is odd, and \emph{white} otherwise.
\item
$\ddd\la\gs2$\\
In this case, we say $\la$ is grey.
\end{description}

From now on we write $\col\la$ for the colour of $\la$. We end this subsection with a \lcnamecref{ddd0nonsplit} we shall need later.

\begin{lemma}\label{ddd0nonsplit}
Suppose $\la\in\pb\tau$ with $\ddd\la=0$, and $\la$ has a $2h$-bar. Then $\lo\la\gs h$.
\end{lemma}

\begin{pf}
Suppose for a contradiction that $\lo\la<h$. Then there are two possibilities.
\begin{description}
\item
$\hi\la=\lo\la+h$\\
In this case the leg lengths of $\la$ are
\[
h-\lo\la+\ope{h-\lo\la}{\lo\la}\tau,\qquad\ope{\lo\la}{\lo\la+h}\tau,
\]
and we claim that these cannot be equal, contradicting the assumption that $\ddd\la=0$. The fact that $\tau$ is an \hbc gives
\[
\ope{h-\lo\la}{\lo\la}\tau\gs\ope{2h-\lo\la}{\lo\la+h}\tau
\]
so in order for the two leg lengths to be equal we would have to have $\ope{\lo\la}{2h-\lo\la}\tau\gs h-\lo\la$. But $\tau$ cannot contain integers $h+i$ and $h-i$ for any $0\ls i<h-\lo\la$, so $\ope{\lo\la}{2h-\lo\la}\tau<h-\lo\la$, a contradiction.
\item
$\hi\la=2h-\lo\la$\\
In this case the leg lengths of $\la$ are 
\[
h-\lo\la+\ope{h-\lo\la}{\lo\la}\tau,\qquad\ope{h-\lo\la}{2h-\lo\la}\tau,
\]
and if these are equal then
\[
h-\lo\la=\ope{\lo\la}{2h-\lo\la}\tau.
\]
But, as in the previous case, this cannot happen.\qedhere
\end{description}
\end{pf}

\subsection{Special partitions}\label{specpartsec}

Our analogue of Richards's formula breaks down for a few (in fact, at most three) \cbvs in $\blk\tau$, and our main theorem will deal with these separately. In order to deal with these vectors, we single out some partitions which we call the \emph{special} partitions in $\pb\tau$. There are up to six of these, though some may be undefined, depending on the value of $\bb\tau$.

\begin{itemize}
\item
If $\bb\tau\ls n-2$, define $\xx\tau$ as follows. Let $a,b\in\{1,\dots,n\}$ be minimal such that $a,b,h-a,h-b\notin\tau$, and define $\xx\tau=\tau\sqcup(h-a,h-b,b,a)$. Note that if $\xx\tau$ is defined, then it is automatically restricted.
\item
If $\bb\tau\ls n-1$, define $\shp\tau$ as follows. Let $a\in\{1,\dots,n\}$ be minimal such that $a,h-a\notin\tau$, and define $\shp\tau=\tau\sqcup(h,h-a,a)$. Note that if $\shp\tau$ is defined, then it is automatically restricted.
\item
Define $\nat\tau$ to be the partition $\tau\sqcup(h,h)$. Observe that $\nat\tau$ is restricted \iff $\tau\neq\varnothing$ (that is, if $\bb\tau\gs1$).
\item
If $\bb\tau\ls n-1$, define $\flt\tau$ as follows. Let $a>h$ be minimal such that $a,2h-a\notin\tau$, and set $\flt\tau=\tau\sqcup(a,2h-a)$.
\item
Define $\ppi\tau$ as follows. Let $a$ be minimal such that $a\in\tau\sqcup(h)$ but $a+h\notin\tau$, and set $\ppi\tau=\tau\sqcup(a+h,h)\setminus(a)$.
\item
If $\bb\tau\gs1$, define $\yy\tau$ as follows. Take the two smallest integers $a<b$ such that $a,b\notin\tau$, $a-h\in\tau$ and $b-h\in\tau\sqcup(a)$. Define $\yy\tau=\tau\sqcup(b,a)\setminus(b-h,a-h)$.
\end{itemize}

It is easy to calculate the values of the function $\ddd$ for these partitions: we get $\ddd\xx\tau=\ddd\shp\tau=\ddd\flt\tau=1$ and $\ddd\nat\tau=\ddd\ppi\tau=\ddd\yy\tau=0$.

\subsection{The main theorem}

Now we can give our main theorem for blocks of \baw $2$.

\begin{thm}\label{mainwt2}
Suppose $\tau$ is an \hbc, and suppose $\la,\mu\in\pb\tau$ with $\mu$ restricted.
\begin{enumerate}[ref=\cref{mainwt2}(\arabic*)]
\item\label{nonspecial}
If $\mu\neq\xx\tau,\shp\tau,\nat\tau$, then there is a partition $\mu^+$ in $\pb\tau$ such that $\mu^+\doms\mu$, $\ddd\mu^+=\ddd\mu$ and $\col\mu=\col\mu^+$. If we let $\mu^+$ be the least dominant such partition, then
\begin{align*}
\dn\la\mu&=
\begin{cases}
1&\text{if }\la=\mu\\
q\text{ or }q^2&\text{if }\mu\domsby\la\domsby\mu^+\text{ and }|\ddd\la-\ddd\mu|=1\\
q^3\text{ or }q^4&\text{if }\la=\mu^+\\
0&\text{otherwise}.
\end{cases}
\\
\intertext{\item
If $\mu=\xx\tau$, then}
\dn\la\mu&=
\begin{cases}
1&\text{if }\la=\xx\tau\\
q&\text{if }\xx\tau\domsby\la\domsby\shp\tau\text{ and }\ddd\la=2\\
q&\text{if }\la=\shp\tau\\
q^2&\text{if }\la=\nat\tau\\
q^3+q&\text{if }\shp\tau\domsby\la\domsby\flt\tau\text{ and }\ddd\la=2\\
q^5+q^3&\text{if }\la=\flt\tau\\
0&\text{otherwise}.
\end{cases}
\\
\intertext{\item
If $\mu=\shp\tau$, then}
\dn\la\mu&=
\begin{cases}
1&\text{if }\la=\shp\tau\\
q&\text{if }\la=\nat\tau\\
q^2&\text{if }\shp\tau\domsby\la\domsby\flt\tau\text{ and }\ddd\la=2\\
q^4+q^2&\text{if }\la=\flt\tau\\
q^2&\text{if }\flt\tau\domsby\la\domsby\ppi\tau\text{ and }\ddd\la=1\\
q^3&\text{if }\la=\ppi\tau\\
0&\text{otherwise}.
\end{cases}
\\
\intertext{\item
If $\mu=\nat\tau$, then}
\dn\la\mu&=
\begin{cases}
1&\text{if }\la=\nat\tau\\
q^3+q&\text{if }\nat\tau\domsby\la\domsby\ppi\tau\text{ and }\ddd\la=1\\
q^2&\text{if }\la=\ppi\tau\\
q^2&\text{if }\ppi\tau\domsby\la\domsby\yy\tau\text{ and }\ddd\la=1\\
q^4&\text{if }\la=\yy\tau\\
0&\text{otherwise}.
\end{cases}
\end{align*}
\end{enumerate}
\end{thm}

We remark that this theorem is enough to determine the $q$-decomposition numbers completely, despite the word ``or'' in part (1): if (for example) $\dn\la\mu$ is given as $q$ or $q^2$, then from \cref{party} we can deduce that $\dn\la\mu=q^2$ \iff $\pry\la=\pry\mu$. In fact it is possible to be more precise: in part (1) of the theorem, if $\dn\la\mu\in\{q,q^2,q^3,q^4\}$, then one can show that $\dn\la\mu\in\{q,q^3\}$ \iff $\la$ has a non-zero part divisible by $h$ but $\mu$ does not.

Most of the remainder of the paper is devoted to proving \cref{mainwt2}. The proof is by induction on $\card\tau$. In \cref{basesec} we deal with the base case, which is where $\tau$ has the form $(l,l-1,\dots,1)$; we treat a large number of \hbcs in the base case in order to avoid having to analyse $[2:1]$-pairs of residue $0$. Our technique here is direct construction of canonical basis vectors starting from known canonical basis vectors. This yields a complete explicit list of all canonical basis vectors when $\tau$ has the form $(l,l-1,\dots,1)$, which we present in the form of two tables.

In \cref{2ksection} we prepare for the inductive step in the proof of \cref{mainwt2}, by taking two blocks $\blk\si$ and $\blk\tau$ forming a $[2:k]$-pair and comparing the dominance order, the functions $\ddd$ and $\col$, and the special partitions in these two blocks. In \cref{indstep} we complete the proof by showing how the canonical bases in $\blk\si$ and $\blk\tau$ differ, and showing that this is compatible with \cref{mainwt2} and the results in \cref{2ksection}.

\section{The base case}\label{basesec}

As with the proof of \cref{mainwt1}, we prove \cref{mainwt2} by induction, taking as base cases the blocks with \bacs of the form $(l,l-1,\dots,1)$.

\begin{basicbox}
\textbf{Assumptions and notation in force for \cref{basesec}}

$\tau=(l,l-1,\dots,1)$ for some $0\ls l\ls n$. Given integers $x_1,\dots,x_s>l$ and $l\gs y_1,\dots,y_t\gs0$, we will write $\rfo{x_1,\dots,x_s}{y_1,\dots,y_t}$ to mean the partition $\tau\sqcup(x_1,\dots,x_s)\setminus(y_1,\dots,y_t)$. 
\end{basicbox}

In the following table, we classify the different types of partition in $\blk\tau$, calculate their $\ddd$-values and give the conditions under which they are restricted.
\needspace{3em}
\[
\begin{array}{cccc}\hline
\la&\text{conditions}&\ddd\la&\text{restricted}\\\hline
\rfo{h-a,h-b,b,a}{\cdot}&l<a<b\ls n&b-a&\text{always}\\
\rfo{h,h}{\cdot}&\text{---}&0&\text{if }l>0\\
\rfo{h+a,h-b,b}{a}&0\ls a\ls l<b\ls n&a+b-l&\text{always}\\
\rfo{h+a,h-a}{\cdot}&1\ls a\ls n&a&\text{always}\\
\rfo{h+b,h+a}{b,a}&0\ls a<b\ls l&b-a-1&\text{if }a\ls l-2\\
\rfo{2h-a,a}{\cdot}&l<a\ls n&a&\text{never}\\
\rfo{2h+a}a&0\ls a\ls l&l-a&\text{never}\\
\hline
\end{array}
\]

Given this, we can find all the partitions in $\pb\tau$ with $\ddd$-value equal to $0$ or $1$, and compute their colours.

\begin{lemma}\label{base0}
The partitions $\la\in\pb\tau$ with $\ddd\la=0$ are given by the following table.
\[
\begin{array}{ccc}\hline
\la&\textup{conditions}&\col\la\\\hline
\rfo{h,h}{\cdot}&\textup{---}&\textup{white}\\
\rfo{h+a+1,h+a}{a+1,a}&0\ls a<l&\left\{\begin{array}l\textup{black if $a$ is even}\\\textup{white if $a$ is odd}\end{array}\right.\\
\rfo{2h+l}l&0\ls l&\left\{\begin{array}l\textup{black if $l$ is even}\\\textup{white if $l$ is odd}\end{array}\right.\\\hline
\end{array}
\]
\end{lemma}

\begin{lemma}\label{base1}
The partitions $\la\in\pb\tau$ with $\ddd\la=1$ are given by the following table.
\[
\begin{array}{ccc}\hline
\la&\textup{conditions}&\col\la\\\hline
\rfo{h-a,h-a-1,a+1,a}{\cdot}&l<a<n&\left\{\begin{array}l\textup{black if $a+l$ is odd}\\\textup{white if $a+l$ is even}\end{array}\right.\\
\rfo{h,h-l-1,l+1}{\cdot}&l<n&\textup{white}\\
\rfo{h+1,h-1}{\cdot}&\text{---}&\textup{grey}\\
\rfo{h+a+2,h+a}{a+2,a}&0\ls a\ls l-2&\textup{grey}\\
\rfo{2h-1,1}{\cdot}&0=l&\textup{grey}\\
\rfo{2h+l-1}{l-1}&1\ls l&\textup{grey}\\\hline
\end{array}
\]
\end{lemma}

Next we need to consider the special partitions in $\pb\tau$. The following result comes directly from the definitions in \cref{specpartsec}.

\needspace{3em}
\begin{lemma}\label{basespec}\indent
\begin{enumerate}
\vspace{-\topsep}
\item
If $l\ls n-2$, the partition $\xx\tau$ equals $\rfo{h-l-1,h-l-2,l+2,l+1}{\cdot}$, and is restricted.
\item
If $l\ls n-1$, the partition $\shp\tau$ equals $\rfo{h,h-l-1,l+1}{\cdot}$, and is restricted.
\item
The partition $\nat\tau$ equals $\rfo{h,h}{\cdot}$, and is restricted \iff $l\gs1$.
\item
If $l\ls n-1$, the partition $\flt\tau$ equals $\rfo{h+1,h-1}{\cdot}$.
\item
The partition $\ppi\tau$ equals $\rfo{h+1,h}1$ if $l\gs1$, or $\rfo{2h}\cdot$ if $l=0$.
\item
If $l\gs2$, the partition $\yy\tau$ equals $\rfo{h+2,h+1}{2,1}$; if $l=1$, then $\yy\tau$ equals $\rfo{2h+1}1$.

\end{enumerate}
\end{lemma}

Now that we have a list of restricted partitions in $\pb\tau$, we can calculate the canonical basis. For each restricted $\mu\in\pb\tau$ we compute $\can\mu$ by applying a suitable combination of operators $f_i$ to a known canonical basis vector $\can\xi$. This partition $\xi$ will have \hbw either $0$ (in which case $\can\xi=\xi$) or $1$ (so that $\can\xi$ is known from \cref{mainwt1}), or will lie in the block $\blk{(l-1,\dots,1)}$ (in which case $\can\xi$ is known by induction on $l$). The results are given in the following table. In each case we give the vector $\can\mu$ in the form of two columns (containing coefficients and partitions), and we give the $\ddd$-value of each partition involved and show how $\can\mu$ is obtained. In each case it is easy to verify that the coefficients $\dn\la\mu$ satisfy \cref{mainwt2}.

We begin with the restricted partitions $\mu\neq\xx\tau,\shp\tau,\nat\tau$.

\begin{longtable}{rllcc}\hline
&$\dn\la\mu$&$\la$&$\ddd\la$&\text{conditions}\\\hline
\endhead
\nextcase
&$1$&\hbox to 190pt{$\rfo{h-l-2,h-l-3,l+3,l+2}\cdot$}&$1$&$l\ls n-3$\\
&$q^2$&$\rfo{h-l-1,h-l-3,l+3,l+1}\cdot$&$2$\\
&$q$&$\rfo{h,h-l-2,l+2}\cdot$&$2$\\
&$q^3$&$\rfo{h,h-l-1,l+1}\cdot$&$1$\\
[3pt]
\rlap{$\can\mu=f_{l+1}\big(\rfo{h-l-2,h-l-3,l+3,l+1}\cdot+q\rfo{h,h-l-2,l+1}\cdot\big)$}\hbox to 15pt{}\\
[3pt]
\hline
\nextcase
&$1$&$\rfo{h-l-1,h-a,a,l+1}\cdot$&$a-l-1$&$l+3\ls a\ls n$\\
&$q^2$&$\rfo{h-l-1,h-a+1,a-1,l+1}\cdot$&$a-l-2$\\
&$q$&$\rfo{h,h-a,a}\cdot$&$a-l$\\
&$q^3$&$\rfo{h,h-a+1,a-1}\cdot$&$a-l-1$\\
[3pt]
\rlap{$\can\mu=f_{a-1}\big(\rfo{h-l-1,h-a,a-1,l+1}\cdot+q\rfo{h,h-a,a-1}\cdot\big)$}\hbox to 15pt{}\\
[3pt]
\hline
\nextcase
&$1$&$\rfo{h-a,h-b,b,a}\cdot$&$b-a$&$l+2\ls a\ls b-2\ls n-2$\\*
&$q^2$&$\rfo{h-a,h-b+1,b-1,a}\cdot$&$b-a-1$\\*
&$q^2$&$\rfo{h-a+1,h-b,b,a-1}\cdot$&$b-a+1$\\*
&$q^4$&$\rfo{h-a+1,h-b+1,b-1,a-1}\cdot$&$b-a$\\[3pt]
\rlap{$\can\mu=f_{a-1}f_{b-1}\rfo{h-a,h-b,b-1,a-1}\cdot$}\hbox to 15pt{}\\[3pt]
\hline
\nextcase
&$1$&$\rfo{h-a,h-a-1,a+1,a}\cdot$&$1$&$l+3\ls a\ls n-1$\\*
&$q^2$&$\rfo{h-a+1,h-a-1,a+1,a-1}\cdot$&$2$\\*
&$q^2$&$\rfo{h-a+2,h-a,a,a-2}\cdot$&$2$\\*
&$q^4$&$\rfo{h-a+2,h-a+1,a-1,a-2}\cdot$&$1$\\[3pt]
\rlap{$\can\mu=f_{a-1}\big(\rfo{h-a,h-a-1,a+1,a-1}\cdot+q^2\rfo{h-a+2,h-a,a-1,a-2}\cdot\big)$}\hbox to 15pt{}\\[3pt]
\hline
\nextcase
&$1$&$\rfo{h+a,h-b,b}a$&$b+a-l$&$0\ls a\ls l-1$\\*
&$q^2$&$\rfo{h+a,h-b+1,b-1}a$&$b+a-l-1$&$l+2\ls b\ls n$\\*
&$q^2$&$\rfo{h+a+1,h-b,b}{a+1}$&$b+a-l+1$\\*
&$q^4$&$\rfo{h+a+1,h-b+1,b-1}{a+1}$&$b+a-l$\\*
[3pt]
\rlap{$\can\mu=f_{b-1}\big(\rfo{h+a,h-b,b-1}{a}+q^2\rfo{h+a+1,h-b,b-1}{a+1}\big)$}\hbox to 15pt{}\\
[3pt]
\hline
\nextcase
&$1$&$\rfo{h+a,h-l-1,l+1}a$&$a+1$&$l\ls n-1$\\*
&$q^2$&$\rfo{h+a,h-a}\cdot$&$a$&$1\ls a\ls l-1$\\*
&$q^2$&$\rfo{h+a+1,h-l-1,l+1}{a+1}$&$a+2$\\*
&$q^4$&$\rfo{h+a+1,h-a-1}\cdot$&$a+1$\\*
[3pt]
\rlap{$\can\mu=f_a\big(\rfo{h+a,h-l-1,l+1}{a+1}+q^2\rfo{h+a,h-a-1}\cdot\big)$}\hbox to 15pt{}\\[3pt]
\hline
\nextcase
&$1$&$\rfo{h+l,h-b,b}l$&$b$&$l+2\ls b\ls n$\\*
&$q^2$&$\rfo{h+l,h-b+1,b-1}l$&$b-1$\\*
&$q^2$&$\rfo{h+b-1,h-b+1}\cdot$&$b-1$\\*
&$q^4$&$\rfo{h+b,h-b}\cdot$&$b$\\[3pt]
\rlap{$\can\mu=f_{b-1}\big(\rfo{h+l,h-b,b-1}l+q^2\rfo{h+b-1,h-b}\cdot\big)$}\hbox to 15pt{}\\[3pt]
\hline
\nextcase
&$1$&$\rfo{h+l,h-l-1,l+1}l$&$l+1$&$1\ls l\ls n-1$\\*
&$q^2$&$\rfo{h+l,h-l}\cdot$&$l$\\*
&$q^4$&$\rfo{h+l+1,h-l-1}\cdot$&$l+1$\\*[3pt]
\rlap{$\can\mu=f_l\rfo{h+l,h-l-1}\cdot$}\hbox to 15pt{}\\[3pt]
\hline
\nextcase
&$1$&$\rfo{h+l,h-l}\cdot$&$l$&$1\ls l\ls n-1$\\
&$q$&$\rfo{h+l,h}l$&$l-1$\\
&$q^2$&$\rfo{h+l+1,h-l-1}\cdot$&$l+1$\\
&$q^2$&$\rfo{2h-l-1,l+1}\cdot$&$l+1$\\
&$q^3$&$\rfo{2h}\cdot$&$l$\\
[3pt]
\rlap{$\can\mu=f_0f_1\dots f_{l-1}\big(\rfo{h{+}l,h{-}l}l+q^2\rfo{h{+}l{+}1,h{-}l{-}1}l+q^2\rfo{2h{-}l{-}1,l{+}1}r+q^4\rfo{2h{-}l}\cdot\big)$}\hbox to 15pt{}\\
[3pt]
\hline
\nextcase
&$1$&$\rfo{h+a,h-a}\cdot$&$a$&$1\ls a\ls l-1$\\*
&$q$&$\rfo{h+a,h}a$&$a-1$\\*
&$q^2$&$\rfo{h+a+1,h-a-1}\cdot$&$a+1$\\*
&$q^3$&$\rfo{h+a+1,h}{a+1}$&$a$\\*
[3pt]
\rlap{$\can\mu=f_a\big(\rfo{h+a,h-a-1}{\cdot}+q\rfo{h+a,h}{a+1}\big)$}\hbox to 15pt{}\\
[3pt]
\hline
\nextcase
&$1$&$\rfo{h+a,h-a}\cdot$&$a$&$l+1\ls a\ls n-1$\\*
&$q^2$&$\rfo{h+a+1,h-a-1}\cdot$&$a+1$\\*
&$q^2$&$\rfo{2h-a-1,a+1}\cdot$&$a+1$\\*
&$q^4$&$\rfo{2h-a,a}\cdot$&$a$\\*
[3pt]
\rlap{$\can\mu=f_a\big(\rfo{h+a,h-a-1}\cdot+q^2\rfo{2h-a-1,a}\cdot\big)$}\hbox to 15pt{}\\
[3pt]
\hline
\nextcase
&$1$&$\rfo{h+n,n+1}\cdot$&$n$&$l\ls n-1$\\*
&$q^4$&$\rfo{h+n+1,n}\cdot$&$n$\\*[3pt]
\rlap{$\can\mu=f_n\rfo{h+n,n}\cdot$}\hbox to 15pt{}\\[3pt]
\hline
\nextcase
&$1$&$\rfo{h+n,n+1}\cdot$&$n$&$l=n$\\*
&$q$&$\rfo{h+n,h}{n}$&$n-1$\\*
&$q^3$&$\rfo{2h}\cdot$&$n$\\*[3pt]
\rlap{$\can\mu=f_0f_1\dots f_n\rfo{h+n}\cdot$}\hbox to 15pt{}\\[3pt]
\hline
\nextcase
&$1$&$\rfo{h+l,h+a}{l,a}$&$l-a-1$&$0\ls a\ls l-2$\\*
&$q^2$&$\rfo{h+l,h+a+1}{l,a+1}$&$l-a-2$\\*
&$q^2$&$\rfo{2h+a}a$&$l-a$\\*
&$q^4$&$\rfo{2h+a+1}{a+1}$&$l-a-1$\\*
[3pt]
\rlap{$\can\mu=f_0^{(2)}\big(\rfo{h+l,h-1}{l,1}+q^2\rfo{2h-1}1\big)$}&&&($a=0$)\\
\rlap{$\can\mu=f_a\big(\rfo{h+l,h+a}{l,a+1}+q^2\rfo{2h+a}{a+1}\big)$}&&&($a\gs1$)\\*
[3pt]
\hline
\nextcase
&$1$&$\rfo{h+l-1,h+l-2}{l-1,l-2}$&$0$&$l\gs2$\\*
&$q^2$&$\rfo{h+l,h+l-2}{l,l-2}$&$1$\\*
&$q^2$&$\rfo{2h+l-1}{l-1}$&$1$\\*
&$q^4$&$\rfo{2h+l}l$&$0$\\
[3pt]
\rlap{$\can\mu=f_{l-1}\big(\rfo{h+l-1,h+l-2}{l,l-2}+q^2\rfo{2h+l-1}l\big)$}\hbox to 15pt{}\\
[3pt]
\hline
\nextcase
&$1$&$\rfo{h+b,h+a}{b,a}$&$b-a-1$&$0\ls a\ls b-2\ls l-3$\\*
&$q^2$&$\rfo{h+b,h+a+1}{b,a+1}$&$b-a-2$\\*
&$q^2$&$\rfo{h+b+1,h+a}{b+1,a}$&$b-a$\\*
&$q^4$&$\rfo{h+b+1,h+a+1}{b+1,a+1}$&$b-a-1$\\*
[3pt]
\rlap{$\can\mu=f_b\big(\rfo{h+b,h+a}{b+1,a}+q^2\rfo{h+b,h+a+1}{b+1,a+1}\big)$}\hbox to 15pt{}\\*
[3pt]
\hline
\nextcase
&$1$&$\rfo{h+a+1,h+a}{a+1,a}$&$0$&$0\ls a\ls l-3$\\*
&$q^2$&$\rfo{h+a+2,h+a}{a+2,a}$&$1$\\*
&$q^2$&$\rfo{h+a+3,h+a+1}{a+3,a+1}$&$1$\\*
&$q^4$&$\rfo{h+a+3,h+a+2}{a+3,a+2}$&$0$\\
[3pt]
\rlap{$\can\mu=f_{a+1}\big(\rfo{h+a+1,h+a}{a+2,a}+q^2\rfo{h+a+3,h+a+1}{a+3,a+2}\big)$}\hbox to 15pt{}\\
[3pt]
\hline
\end{longtable}

Now we do the same for $\can{\xx\tau}$, $\can{\shp\tau}$ and $\can{\nat\tau}$. In each case we identify the special partitions appearing, and we can complete the verification of \cref{mainwt2} in the base cases.

\begin{longtable}{rlllcc}\hline
&$\dn\la\mu$&$\la$&&$\ddd\la$&\text{conditions}\\\hline
\endhead
\nextcase
&$1$&$\rfo{h-l-1,h-l-2,l+2,l+1}{\cdot}$&$\xx\tau$&$1$&$1\ls l\ls n-2$\\*
&$q$&$\rfo{h,h-l-2,l+2}{\cdot}$&&$2$\\*
&$q$&$\rfo{h,h-l-1,l+1}{\cdot}$&$\shp\tau$&$1$\\*
&$q^2$&$\rfo{h,h}{\cdot}$&$\nat\tau$&$0$\\*
&$q^3+q$&$\rfo{h+1,h-l-1,l+1}1$&&$2$\\*
&$q^5+q^3$&$\rfo{h+1,h-1}{\cdot}$&$\flt\tau$&$1$\\*
[3pt]
\rlap{$\can\mu=f_0\dots f_l\big(\rfo{h-l-1,h-l-2,l+2}{\cdot}+q\rfo{h,h-l-1}{\cdot}\big)$}\hbox to 15pt{}\\
[3pt]
\hline
\nextcase
&$1$&$\rfo{h-1,h-2,2,1}{\cdot}$&$\xx\tau$&$1$&$l=0\ls n-2$\\*
&$q$&$\rfo{h,h-2,2}{\cdot}$&&$2$\\*
&$q$&$\rfo{h,h-1,1}{\cdot}$&$\shp\tau$&$1$\\*
&$q^2$&$\rfo{h,h}{\cdot}$&$\nat\tau$&$0$\\*
&$q^5+q^3$&$\rfo{h+1,h-1}{\cdot}$&$\flt\tau$&$1$\\*
[3pt]
\rlap{$\can\mu=f_0\big(\rfo{h-1,h-2,2}{\cdot}+q\rfo{h,h-1}{\cdot}\big)$}\hbox to 15pt{}\\
[3pt]
\hline
\nextcase
&$1$&$\rfo{h,h-l-1,l+1}\cdot$&$\shp\tau$&$1$&$1\ls l\ls n-1$\\*
&$q$&$\rfo{h,h}{\cdot}$&$\nat\tau$&$0$\\*
&$q^2$&$\rfo{h+1,h-l-1,l+1}1$&&$2$\\*
&$q^4+q^2$&$\rfo{h+1,h-1}{\cdot}$&$\flt\tau$&$1$\\*
&$q^3$&$\rfo{h+1,h}1$&$\ppi\tau$&$0$\\*
[3pt]
\rlap{$\can\mu=f_0^{(2)}\big(\rfo{h-1,h-l-1,l+1}1+q\rfo{h,h-1}1\big)$}\hbox to 15pt{}\\*
[3pt]
\hline
\nextcase
&$1$&$\rfo{h,h-1,1}{\cdot}$&$\shp\tau$&$1$&$l=0$\\*
&$q$&$\rfo{h,h}{\cdot}$&$\nat\tau$&$0$\\*
&$q^4+q^2$&$\rfo{h+1,h-1}{\cdot}$&$\flt\tau$&$1$\\*
&$q^2$&$\rfo{2h-1,1}{\cdot}$&&$1$\\*
&$q^3$&$\rfo{2h}{\cdot}$&$\ppi\tau$&$0$\\*
[3pt]
\rlap{$\can\mu=f_0\big(\rfo{h,h-1}{\cdot}+q^2\rfo{2h-1}{\cdot}\big)$}\hbox to 15pt{}\\
[3pt]
\hline
\nextcase
&$1$&\hbox to 190pt{$\rfo{h,h}{\cdot}$}&$\nat\tau$&$0$&$l\gs2$\\*
&$q^3+q$&$\rfo{h+1,h-1}{\cdot}$&&$1$\\*
&$q^2$&$\rfo{h+1,h}1$&$\ppi\tau$&$0$\\*
&$q^2$&$\rfo{h+2,h}2$&&$1$\\*
&$q^4$&$\rfo{h+2,h+1}{2,1}$&$\yy\tau$&$0$\\*
[3pt]
\rlap{$\can\mu=f_0^{(2)}\big(\rfo{h,h-1}1+q^2\rfo{h+2,h-1}{2,1}\big)$}\hbox to 15pt{}\\*
[3pt]
\hline
\nextcase
&$1$&$\rfo{h,h}{\cdot}$&$\nat\tau$&$0$&$l=1$\\*
&$q^3+q$&$\rfo{h+1,h-1}{\cdot}$&&$1$\\*
&$q^2$&$\rfo{h+1,h}1$&$\ppi\tau$&$0$\\*
&$q^2$&$\rfo{2h}{\cdot}$&&$1$\\*
&$q^4$&$\rfo{2h+1}1$&$\yy\tau$&$0$\\*
[3pt]
\rlap{$\can\mu=f_0^{(2)}\big(\rfo{h,h-1}1+q^2\rfo{2h-1}1\big)$}\hbox to 15pt{}\\*
[3pt]
\hline
\end{longtable}

\section{$[2:k]$-pairs}\label{2ksection}

Now we come to the inductive step in the proof of \cref{mainwt2}. For this, we need to study $[2:k]$-pairs in more detail.

Suppose $\blk\si$ and $\blk\tau$ form a $[2:k]$-pair of residue $i$. The results in \cref{basesec} allow us to avoid $[2:1]$-pairs of residue $0$, which are difficult to work with, so we concentrate on the other cases.

\begin{basicbox}
\textbf{Assumptions and notation in force for \cref{2ksection}}

$\sigma$ and $\tau$ are \hbcs, and $\blk\sigma$ and $\blk\tau$ form a $[2:k]$-pair of residue $i$, where either $i\gs1$ or $k\gs3$.
\end{basicbox}

In this section we compare the dominance order, the functions $\ddd$ and $\col$ and the special partitions in $\blk\si$ and $\blk\tau$. The most difficult cases are where $\blk\si$ and $\blk\tau$ are not Scopes--Kessar equivalent. From \cref{skequiv}, $\blk\si$ and $\blk\tau$ are Scopes--Kessar equivalent \iff one of the following holds:
\begin{itemize}
\item
$k\gs5$;
\item
$1\ls i<n$ and $k\gs2$;
\item
$i=n$.
\end{itemize}

\subsection{Exceptional partitions}\label{21sec}

First we look at the case where $\blk\si$ and $\blk\tau$ are not Scopes--Kessar equivalent, and study the exceptional partitions for the pairs $(\blk\si,\blk\tau)$.

We begin with the case where $1\ls i<n$ and $k=1$. Consider runners $-i-1,-i,i,i+1$ in the \abd for $\sigma$; let $b_{-i-1},b_{-i},b_i,b_{i+1}$ be the positions of the lowest beads on each of these runners. Then (because $\sigma$ is an \hbc) $b_{-i-1}+b_{i+1}=b_{-i}+b_i=-h$, and (because $\sigma$ has exactly one addable $i$-node) $b_i=b_{i+1}+h-1$. Now we consider three possibilities, depending on the relative order of $b_{-i-1},b_{-i},b_i,b_{i+1}$:
\begin{description}
\item[Type A]
The addable $i$-node of $\sigma$ lies in column $ah+i+1$ for some $a\gs1$, in which case $b_{-i}<b_{-i-1}<b_{i+1}<b_i$;
\item[Type B]
The addable $i$-node of $\sigma$ lies in column $i+1$, in which case $b_{-i}<b_{i+1}<b_{-i-1}<b_i$;
\item[Type C]
The addable $i$-node of $\sigma$ lies in column $ah-i$ for some $a\gs1$, in which case $b_{i+1}<b_i<b_{-i}<b_{-i-1}$.
\end{description}

For example, take $h=7$. We illustrate $[2:1]$-pairs of the three possible types, with $i=1$ in each case.

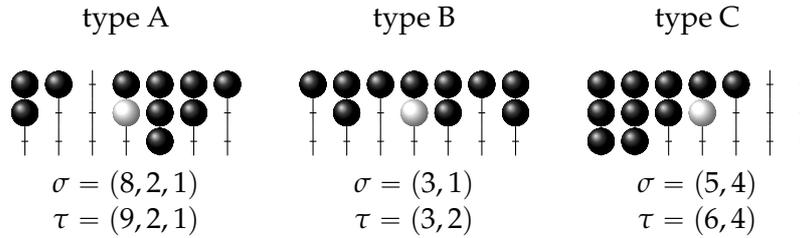
\newcommand\minusone{-1}
\newcommand\minustwo{-2}
\newcommand\minusthree{-3}
\[
\begin{array}{c@{\qquad}c@{\qquad}c}
\text{type A}&\text{type B}&\text{type C}\\[12pt]
\abacus(bbnbbbb,bnnxbbn,nnnnbnn)
&
\abacus(bbbbbbb,nbnxbnb,nnnnnnn)
&
\abacus(bbbbbnn,bbbxnnn,bbnnnnn)\\[12pt]
\sigma=(8,2,1)&\sigma=(3,1)&\sigma=(5,4)\\
\tau=(9,2,1)&\tau=(3,2)&\tau=(6,4)
\end{array}
\]

Now we can classify and study the exceptional partitions in $\blk\sigma$ and $\blk\tau$. We use the abacus notation for partitions introduced in \cref{absec}.

\begin{propn}\label{21stuff}
Suppose $\blk\sigma$ and $\blk\tau$ form a $[2:1]$-pair of residue $i$, where $1\ls i<n$.
\begin{enumerate}[ref=(\arabic*)]
\item\label{21stuffexcnu}
There are three exceptional partitions $\pb\sigma$, which can be written $\alpha\domsby\beta\domsby\gamma$, and there are three exceptional partitions in $\pb\tau$, which can be written $\halpha\domsby\hbeta\domsby\hgamma$. These partitions have abacus notation given by the following table.
\[
\begin{array}{lccc@{\qquad}ccc}\hline
&\alpha&\beta&\gamma&\halpha&\hbeta&\hgamma\\\hline
\text{Type A}&\ba{-i}&\ba{i,i+1}&\ba{i+1}&\ba{i}&\ba{i,i+1}&\ba{-i-1}\\
\text{Type B}&\ba{i,i+1}&\ba{-i}&\ba{i+1}&\ba{i}&\ba{-i-1}&\ba{i,i+1}\\
\text{Type C}&\ba{i+1}&\ba{i,i+1}&\ba{-i}&\ba{-i-1}&\ba{i,i+1}&\ba{i}\\\hline
\end{array}
\]
\Item\label{21stuffpsi}
\begin{alignat*}3
\psi_i(\alpha)&=\halpha,&\qquad\psi_i(\beta)&=\hgamma,&\qquad\psi_i(\gamma)&=\hbeta,
\\
f_i\alpha&=q^{-2}\halpha+\hbeta,&\qquad
f_i\beta&=\halpha+\hgamma&\qquad
f_i\gamma&=\hbeta+q^2\hgamma.
\end{alignat*}
\item\label{sigma1cbv}
$\alpha$ and $\halpha$ are restricted, with
\[
\can\alpha=\alpha+q^2\beta+q^4\gamma,\qquad\can\halpha=\halpha+q^2\hbeta+q^4\hgamma.
\]
\end{enumerate}
\end{propn}

\needspace{5em}
\begin{pfenum}
\item
By examining runners $i,i+1$ of the abacus, we find that the exceptional partitions in $\pb\sigma$ are precisely $\ba{-i}$, $\ba{i,i+1}$ and $\ba{i+1}$, and similarly for $\pb\tau$. The following diagrams show runners $i$ and $i+1$ of the abacus.
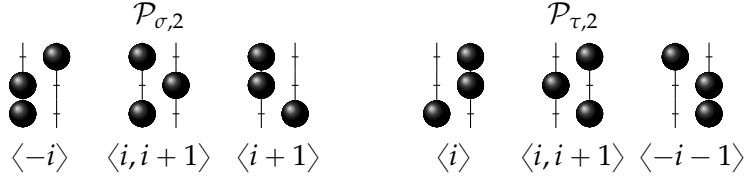
\[
\begin{array}{c@{\qquad\quad}c}
\pb\sigma&\pb\tau\\[3pt]
\begin{array}{c@{\qquad}c@{\qquad}c}
\abacus(nb,bn,bn)
&
\abacus(bn,nb,bn)
&
\abacus(bn,bn,nb)\\[12pt]
\mathclap{\ba{-i}}&\mathclap{\ba{i,i+1}}&\mathclap{\ba{i+1}}
\end{array}
&
\begin{array}{c@{\qquad}c@{\qquad}c}
\abacus(nb,nb,bn)
&
\abacus(nb,bn,nb)
&
\abacus(bn,nb,nb)\\[12pt]
\mathclap{\ba{i}}&\mathclap{\ba{i,i+1}}&\mathclap{\ba{-i-1}}
\end{array}
\end{array}
\]
The dominance ordering on these triples of partitions is easily checked from the abacus.
\item
These statements follow by considering the arrangement of addable and removable $i$-nodes for $\alpha,\beta,\gamma$. For example, $\gamma$ has two addable $i$-nodes and one removable $i$-node, with the removable $i$-node to the right of the addable $i$-nodes. Adding the leftmost addable $i$-node yields $\hbeta$, while adding the other addable $i$-node yields $\hgamma$. Hence $\psi_i(\gamma)=\hbeta$, and $f_i\gamma=\hbeta+q^2\hgamma$.
\item
Let $\delta$ be the \hsp obtained by removing the unique removable $i$-node from any of $\alpha$, $\beta$ or $\gamma$. Then $\delta$ is an \hbc, so $\can\delta=\delta$. So the vector
\[
f_i\delta=\alpha+q^2\beta+q^4\gamma
\]
is bar-invariant, and therefore equals $\can\alpha$ (and $\alpha$ is necessarily restricted). Now from part \ref{21stuffpsi} we can calculate
\[
f_i\can\alpha=(q^{-2}+q^2)(\halpha+q^2\hbeta+q^4\hgamma).
\]
Hence the vector $\halpha+q^2\hbeta+q^4\hgamma$ is bar-invariant, and so must equal $\can\halpha$.
\end{pfenum}

Now we give the corresponding results for the case $i=0$ and $k=3$.

\needspace{3em}
\begin{propn}\label{23stuff}
Suppose $\blk\sigma$ and $\blk\tau$ form a $[2:3]$-pair of residue $0$.
\begin{enumerate}[ref=(\arabic*)]
\item\label{23stuffexcnu}
The exceptional partitions in $\pb\sigma$ can be written $\alpha\domsby\beta\domsby\gamma$, where
\begin{alignat*}3
\alpha&=\ba1,&\qquad\beta&=\ba{0,1},&\qquad\gamma&=\ba0. 
\\
\intertext{The exceptional partitions in $\pb\tau$ can be written $\halpha\domsby\hbeta\domsby\hgamma$, where}
\halpha&=\ba0,&\qquad\hbeta&=\ba{0,1},&\qquad\hgamma&=\ba{-1}. 
\end{alignat*}
\Item\label{23stuffpsi}
\[
\psi_0(\alpha)=\halpha,\qquad\psi_0(\beta)=\hgamma,\qquad\psi_0(\gamma)=\hbeta,
\]
\begin{align*}
f_0^{(3)}\alpha&=q^{-3}\halpha+q^{-1}\hbeta\\
f_0^{(3)}\beta&=(1+q^{-2})\halpha+\hbeta+(q^2+1)\hgamma\\
f_0^{(3)}\gamma&=\halpha+(q^2+1)\hbeta+(q^4+q^2)\hgamma.
\end{align*}
\item\label{gminusi}
$\alpha$ and $\halpha$ are restricted, with
\[
\can\alpha=\alpha+q\beta+q^3\gamma,\qquad\can\halpha=\halpha+q^2\hbeta+q^4\hgamma.
\]
\end{enumerate}
\end{propn}

\needspace{5em}
\begin{pfenum}
\item
This is a matter of checking possible abacus configurations. Runners $-1,0,1$ of the \abds for the exceptional partitions are shown below.
\[
\begin{array}{c@{\qquad\quad}c}
\pb\sigma&\pb\tau\\[3pt]
\begin{array}{c@{\qquad}c@{\qquad}c}
\abacus(bbn,bbn,nxb,bnn,bnn)
&
\abacus(bbn,bnb,bxn,nbn,bnn)
&
\abacus(bnb,bbn,bxn,bnn,nbn)\\
\mathclap{\ba1}&\mathclap{\ba{0,1}}&\mathclap{\ba0}
\end{array}
&
\begin{array}{c@{\qquad}c@{\qquad}c}
\abacus(bnb,nbb,nxb,nnb,nbn)
&
\abacus(nbb,bnb,nxb,nbn,nnb)
&
\abacus(nbb,nbb,bxn,nnb,nnb)\\
\mathclap{\ba0}&\mathclap{\ba{0,1}}&\mathclap{\ba{-1}}
\end{array}
\end{array}
\]
\item
As in \cref{21stuff}, these statements follow by considering the configuration of addable and removable $0$-nodes for each of $\alpha$, $\beta$ and $\gamma$.

As an example, we show how to compute the coefficient of $\halpha$ in $f_0^{(3)}\beta$. As we can see from the abacus display, the addable and removable $0$-nodes of $\beta=\ba{0,1}$ consist of (from left to right):
\begin{itemize}
\item
an addable node in column $1$;
\item
a removable node in column $h$;
\item
addable nodes in columns $h+1$, $2h$ and $2h+1$.
\end{itemize}
There are three possible ways to add three of the addable $0$-nodes, yielding the partitions $\halpha$, $\hbeta$ and $\hgamma$. In particular, $\halpha$ is obtained by adding the addable nodes in columns $1$, $h+1$ and $2h$. Since $h$ occurs exactly once as a part of $\beta$, the definition of the action of $f_0^{(3)}$ gives a coefficient of $q^{-2}\times(1+q^2)=(1+q^{-2})$.
\item
This is proved as in \cref{21stuff}, using the \hbc $\delta=\sigma\sqcup(2h-1)$.
\end{pfenum}

\subsection{$\ddd$-values and dominance}\label{2ksec}

Now we consider how the dominance orders and the $\ddd$-values differ for corresponding \hsps in two blocks forming a $[2:k]$-pair.

First we consider $\ddd$-values and colours. We start with unexceptional partitions.

\begin{propn}\label{unexddd}
Suppose $\la\in\pb\sigma$ is unexceptional. Then $\ddd\psi_i(\la)=\ddd\la$.
\end{propn}

\begin{pf}
For this proof, we write $\mu=\psi_i(\la)$. To prove that $\ddd\la=\ddd\mu$, we prove the stronger statement that $\mu$ has the same leg lengths as $\la$. Let $\lo\la<\hi\la$ be the bar-positions for $\la$. Define $\bea\la$ to be the set of occupied positions in the \abd for $\la$, and set
\[
\hbea\la=\bea\la\setminus\{\pm \lo\la,\pm(\lo\la-h),\pm \hi\la,\pm(\hi\la-h)\}\cup\{0\}.
\]
define $\lo\mu$, $\hi\mu$ and $\hbea\mu$ similarly.

Now observe that the leg lengths of $\la$ are simply
\[
\card{\hbea\la\cap(\lo\la-h,\lo\la)},\qquad\card{\hbea\la\cap(\hi\la-h,\hi\la)}
\]
(where we use the usual notation $(x,y)$ for the open interval between $x$ and $y$). The leg lengths of $\mu$ can be expressed similarly, so we need to compare the sets $\hbea\la\cap(\lo\la-h,\lo\la)$ and $\hbea\mu\cap(\lo\mu-h,\lo\mu)$ (and the corresponding sets for $\hi\la$ and $\hi\mu$).

We consider the cases $i\neq0$ and $i=0$ separately. First suppose $i\neq0$. In this case, we let $\phi_i$ denote the bijection from $\bbz$ to $\bbz$ given by
\[
\phi_i(c)=
\begin{cases}
c+1&\text{if }c\equiv i\text{ or }-i-1\ppmod h\\
c-1&\text{if }c\equiv i+1\text{ or }-i\ppmod h\\
c&\text{otherwise}.
\end{cases}
\]

Then $\lo\mu=\phi_i(\lo\la)$ unless either
\begin{enumerate}[label=(\alph*),ref=(\alph*)]
\item\label{aln}
$\lo\la=n+1$ and $i=n$, or
\item\label{aabb}
$\lo\la=\hi\la-1\equiv i$ or $-i-1\ppmod h$.
\end{enumerate}
In either of these cases $\lo\mu=\lo\la$. Similarly $\hi\mu=\phi_i(\hi\la)$ except in case \ref{aabb} above, in which case $\hi\mu=\hi\la$.

In any case, $\hbea\mu=\phi_i(\hbea\la)$; so since $\phi_i$ is a bijection on $\bbz$, to compare the sizes of the sets $\hbea\la\cap(\lo\la-h,\lo\la)$ and $\hbea\mu\cap(\lo\mu-h,\lo\mu)$ we just need to consider the possible integers $c$ such that $c\in\hbea\la\cap(\lo\la-h,\lo\la)$ and $\phi_i(c)\notin(\lo\mu-h,\lo\mu)$, or $c\in\hbea\mu\cap\notin(\lo\mu-h,\lo\mu)$ and $\phi_i(c)\notin(\lo\la-h,\lo\la)$.

Consider first the case $\lo\mu=\lo\la$. The only way an integer $c$ can satisfy $c\in(\lo\la-h,\lo\la)\not\ni\phi_i(c)$ is if
\begin{itemize}
\item
$c=\lo\la-1$ and $\lo\la\equiv i+1$ or $-i\ppmod h$, or
\item
$c=\lo\la-h+1$ and $\lo\la\equiv i$ or $-i-1\ppmod h$.
\end{itemize}
But if $\lo\la\equiv i+1$ or $-i\ppmod h$ then $\phi_i(\lo\la)=\lo\la-1$, so the only way we can have $\lo\mu=\lo\la$ is in case \ref{aln} above; that is, $\lo\la=n+1=i+1$. But in that case $\lo\la-1=h-\lo\la\notin\hbea\la,\hbea\mu$.

Similarly if $\lo\la\equiv i$ or $-i-1\ppmod h$ then $\phi_i(\lo\la)=\lo\la+1$, so the only way we can have $\lo\mu=\lo\la$ is in case \ref{aabb} above. But then $\lo\la-h+1=\hi\la-h\notin\hbea\la,\hbea\mu$.

So $\phi_i$ maps $\hbea\la\cap(\lo\la-h,\lo\la)$ bijectively to $\hbea\mu\cap(\lo\mu-h,\lo\mu)$.

Now consider the case $\lo\mu\neq \lo\la$. This can happen in either of two ways.
\begin{description}
\item
$\lo\la\equiv i$ or $-i-1\ppmod h$ and $\lo\mu=\lo\la+1$\\
Observe that if $c\in(\lo\la-h,\lo\la)$, then $\phi_i(c)\in(\lo\la-h+1,\lo\la+1)$ unless $c=\lo\la-h+1$. Similarly, if $\phi_i(c)\in(\lo\la-h+1,\lo\la+1)$, then $c\in(\lo\la-h,\lo\la)$ unless $c=\lo\la$. So to show that $\card{\hbea\la\cap(\lo\la-h,\lo\la)}=\card{\hbea\mu\cap(\lo\mu-h,\lo\mu)}$, we just need to show that $\lo\la-h+1\in\hbea\la$ \iff $\lo\la\in\hbea\mu$; that is, $\lo\la-h+1\in\hbea\la$ \iff $\lo\la+1\in\hbea\la$. If $\lo\la-h+1\notin\hbea\la\ni \lo\la+1$, then $\hi\la=\lo\la+1$, so that $\lo\mu=\lo\la$, contrary to assumption. On the other hand, if $\lo\la-h+1\in\hbea\la\not\ni \lo\la+1$, then one of the following must occur on runners $i,i+1$ of the \abd for $\la$.
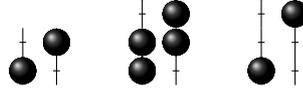
\[
\abacus(nb,bn)\qquad\abacus(nb,bb,bn)\qquad\abacus(nb,nn,bn)
\]
But then $\la$ is exceptional, contrary to assumption.
\item
$\lo\la\equiv i+1$ or $-i\ppmod h$ and $\lo\mu=\lo\la-1$
\\
This case is similar to the previous one; here we must show that $\lo\la-h-1\in\hbea\la$ \iff $\lo\la-1\in\hbea\la$, and this is done in a similar way.
\end{description}
So in all cases, $\card{\hbea\la\cap(\lo\la-h,\lo\la)}=\card{\hbea\mu\cap(\lo\mu-h,\lo\mu)}$.

\smallskip
Now we consider the case $i=0$, where we take the same approach. We define the bijection $\phi_0:\bbz\to\bbz$ by
\[
\phi_0(c)=
\begin{cases}
c+2&\text{if }c\equiv-1\ppmod h\\
c-2&\text{if }c\equiv1\ppmod h\\
c&\text{otherwise}.
\end{cases}
\]

Suppose first that $\lo\la\nequiv\pm1\ppmod h$. Then $\lo\mu=\lo\la$, and the only way we can have an integer $c$ with $c\in(\lo\la-h,\lo\la)\not\ni\phi_0(c)$ is if $\lo\la\equiv0\ppmod h$ and $c=\lo\la-1$ or $\lo\la-h+1$. So if $\lo\la\nequiv0\ppmod h$, then $\phi_0$ maps $\hbea\la\cap(\lo\la-h,\lo\la)$ bijectively to $\hbea\mu\cap(\lo\la-h,\lo\la)$. If $\lo\la\equiv0\ppmod h$, then the assumption that $\la$ is unexceptional means that $\lo\la-1,\lo\la-h-1\in\hbea\la$ while $\lo\la+1,\lo\la-h+1\notin\hbea\la$. So
\begin{align*}
\hbea\la\cap(\lo\la-h,\lo\la)&=\hbea\la\cap(\lo\la-h+1,\lo\la-1)\cup\{\lo\la-1\}\\
\hbea\mu\cap(\lo\la-h,\lo\la)&=\hbea\la\cap(\lo\la-h+1,\lo\la-1)\cup\{\lo\la-h+1\}
\end{align*}
and hence $\card{\hbea\la\cap(\lo\la-h,\lo\la)}=\card{\hbea\mu\cap(\lo\la-h,\lo\la)}$.

Now suppose $\lo\la\equiv-1\ppmod h$. Then the assumption that $\la$ is unexceptional means that $\hi\la$ does not equal $\lo\la+1$ or $\lo\la+2$, so $\lo\mu=\phi_0(\lo\la)$. If $c\in(\lo\la-h,\lo\la)$, then $\phi_0(c)\in(\lo\la-h+2,\lo\la+2)$ unless $c=\lo\la-h+1$ or $\lo\la-h+2$, while if $c\in(\lo\la-h+2,\lo\la+2)$ then $\phi_0(c)\in(\lo\la-h,\lo\la)$ unless $c=\lo\la$ or $\lo\la+2$. So in order to show $\card{\hbea\la\cap(\lo\la-h,\lo\la)}=\card{\hbea\mu\cap(\lo\la-h,\lo\la)}$, we must show that
\begin{align*}
\card{\hbea\la\cap\{\lo\la-h+1,\lo\la-h+2\}}&=\card{\hbea\mu\cap\{\lo\la,\lo\la+1\}},\\
\intertext{that is}
\card{\hbea\la\cap\{\lo\la-h+1,\lo\la-h+2\}}&=\card{\hbea\la\cap\{\lo\la+1,\lo\la+2\}}.
\end{align*}
This follows from the fact that $\la$ is unexceptional by considering possible \abds.

\smallskip
So in all cases (for $i\neq0$ and $i=0$) $\card{\hbea\la\cap(\lo\la-h,\lo\la)}=\card{\hbea\mu\cap(\lo\mu-h,\lo\mu)}$. In the same way we can show that $\card{\hbea\la\cap(\hi\la-h,\hi\la)}=\card{\hbea\mu\cap(\hi\mu-h,\hi\mu)}$.
\end{pf}

Now we do the same for colours.

\begin{propn}\label{unexcl}
Suppose $\la\in\pb\sigma$ is unexceptional. Then $\col\psi_i(\la)=\col\la$.
\end{propn}

\begin{pf}
Again, we write $\mu=\psi_i(\la)$. Recall that $\bb\sigma$ denotes the number of parts of $\sigma$ which are less than $h$. Observe that $\bb\sigma=\bb\tau$: for $i\gs1$ this is completely clear, while if $i=0$ then the assumption that $k\gs3$ means that $h-1\in\sigma\not\ni1$, and the set $\tau\cap\{1,\dots,h-1\}$ is obtained from $\sigma\cap\{1,\dots,h-1\}$ by replacing $h-1$ with $1$.

By \cref{unexddd} $\ddd\la=\ddd\mu$, so if $\ddd\la\gs2$ the result is trivial. We consider the two remaining cases.

\begin{description}
\item
$\ddd\la=0$
\\
Let $l$ be the (repeated) leg length of $\la$ (and of $\mu$, from the proof of \cref{unexddd}).

Suppose $\la$ has a $2h$-bar. By \cref{ddd0nonsplit} this means that $h\ls\lo\la=\hi\la-h$. The definition of colour means that $\la$ is black \iff either
\begin{itemize}
\item
$\lo\la\notin\sigma$ and $l+\bb\sigma$ is even, or
\item
$\lo\la\in\sigma$ and $l+\bb\sigma$ is odd.
\end{itemize}
Now $\phi_i(\lo\la)$ and $\phi_i(\hi\la)$ are the bar positions of $\mu$, and $\phi_i(\lo\la)\in\tau$ \iff $\lo\la\in\sigma$. Now the fact that $\bb\tau=\bb\sigma$ means that $\la$ is black \iff $\mu$ is.

Alternatively, suppose $\la$ has two $h$-bars. Then $\la$ is black \iff $l+\bb\sigma$ is odd. Clearly also $\mu$ has two $h$-bars, and so $\mu$ is black \iff $l+\bb\tau=l+\bb\sigma$ is odd.
\item
$\ddd\la=1$
\\
In this case let $l,l+1$ be the leg lengths of $\la$ and of $\mu$.

If $\hi\la>h$, then we claim that $\hi\mu>h$: the only way this could potentially fail is if $i=0$ and $\hi\la=h+1$; but then $\la$ would be exceptional. So $\la$ and $\mu$ are both grey if $\hi\la>h$.

Similarly, if $\hi\la\ls h$ then $\hi\mu\ls h$. In this case $\la$ and $\mu$ are both black if $l+\bb\sigma=l+\bb\tau$ is odd, and both white otherwise.\qedhere
\end{description}
\end{pf}

Next we compare the dominance orders in $\pb\sigma$ and $\pb\tau$.

\begin{propn}\label{unexdom}
Suppose $\la,\mu\in\pb\sigma$ are both unexceptional and $|\ddd\la-\ddd\mu|\ls1$. Then $\la\domby\mu$ \iff $\psi_i(\la)\domby\psi_i(\mu)$.
\end{propn}

\begin{pf}
By \cref{domddd} $\la$ and $\mu$ are comparable in the dominance order. By \cref{unexddd} $|\ddd\psi_i(\la)-\ddd\psi_i(\mu)|\ls1$, so $\psi_i(\la)$ and $\psi_i(\mu)$ are also comparable in the dominance order. Hence
\[
\la\domby\mu\ \Longleftrightarrow\ \la\lslex\mu\ \Longleftrightarrow\ \psi_i(\la)\lslex\psi_i(\mu)\ \Longleftrightarrow\ \psi_i(\la)\domby\psi_i(\mu),
\]
with the middle implication following from \cref{unexlex}\ref{unexlexlex}.
\end{pf}

\begin{rmk}
In fact, \cref{unexdom} is true even without the hypothesis that $|\ddd\la-\ddd\mu|\ls1$, but this is harder to prove, and we do not need it.
\end{rmk}

Now we come to the exceptional partitions in the cases where $\blk\si$ and $\blk\tau$ are not Scopes--Kessar equivalent. In these cases we use the labelling for the exceptional partitions introduced in \cref{21stuff,23stuff}.

\begin{propn}\label{domexc}
Suppose that either $i=0$ and $k=3$, or $1\ls i<n$ and $k=1$. Let $\alpha\domsby\beta\domsby\gamma$ be the exceptional partitions in $\pb\sigma$, and $\halpha\domsby\hbeta\domsby\hgamma$ the exceptional partitions in $\pb\tau$. Then there is an integer $d\gs1$ such that the following hold.
\begin{enumerate}[ref=(\arabic*)]
\item\label{da}
$\ddd\alpha=\ddd\hbeta=\ddd\gamma=d$ and $\col\alpha=\col\hbeta=\col\gamma$.
\item\label{dha}
$\ddd\halpha=\ddd\beta=\ddd\hgamma=d-1$ and $\col\halpha=\col\beta=\col\hgamma$.
\item\label{aghb}
If $\la\in\pb\sigma$ is unexceptional with $|\ddd\la-d|\ls1$, then either
\begin{itemize}
\item
$\la\doms\gamma$ and $\psi_i(\la)\doms\hbeta$, or
\item
$\la\domsby\alpha$ and $\psi_i(\la)\domsby\hbeta$.
\end{itemize}
\item\label{hahgb}
If $\la\in\pb\sigma$ is unexceptional with $|\ddd\la-(d-1)|\ls1$, then either
\begin{itemize}
\item
$\la\doms\beta$ and $\psi_i(\la)\doms\hgamma$, or
\item
$\la\domsby\beta$ and $\psi_i(\la)\domsby\halpha$.
\end{itemize}
\end{enumerate}
\end{propn}

\begin{pf}
To prove parts \ref{da} and \ref{dha}, we define two integers $l$ and $m$ depending on the abacus configuration of $\sigma$. We continue to use the notation $\ope xy\sigma$ for the number of parts of $\sigma$ lying strictly between~$x$ and~$y$. We identify five cases.
\begin{enumerate}[label=(\alph*)]
\item\label{low}
Suppose $i\gs1$ and the addable $i$-node of $\sigma$ lies in column $i+1$. Then define
\begin{alignat*}2
l&=i+\ope{i+1}{h-i-1}\sigma,&\qquad m&=\ope{i+1}{h+i}\sigma.
\\
\intertext{
\item
Suppose $i\gs1$ and the addable $i$-node of $\sigma$ lies in column $h-i$. Then define}
l&=i+\ope{i+1}{h-i-1}\sigma,&\qquad m&=\ope{h-i}{2h-i-1}\sigma.
\\
\intertext{
\item
Suppose $i\gs1$ and the addable $i$-node of $\sigma$ lies in column $ah+i+1$ for $a\gs1$. Then define}
l&=\ope{(a-1)h+i+1}{ah+i-1}\sigma,&\qquad m&=\ope{ah+i+1}{(a+1)h+i-1}\sigma.
\\
\intertext{
\item
Suppose $i\gs1$ and the addable $i$-node of $\sigma$ lies in column $ah+h-i$ for $a\gs2$. Then define}
l&=\ope{ah-i}{(a+1)h-i-2}\sigma,&\qquad m&=\ope{(a+1)h-i}{(a+2)h-i-2}\sigma.
\\
\intertext{
\item
Suppose $i=0$. Define}
l&=\ope1{h-1}\sigma,&\qquad m&=\ope{h+1}{2h-1}\sigma.
\end{alignat*}
\end{enumerate}

Now by checking the partitions $\alpha,\beta,\gamma,\halpha,\hbeta,\hgamma$ in each of the five cases, we find that
\begin{alignat*}2
\alpha&\text{ has leg lengths }l+1,m&\qquad\halpha&\text{ has leg lengths }l+1,m+1\\
\beta&\text{ has leg lengths }l,m&\qquad\hbeta&\text{ has leg lengths }l+1,m\\
\gamma&\text{ has leg lengths }l+1,m&\qquad\hgamma&\text{ has leg lengths }l,m,
\end{alignat*}
so that, taking $d=l-m+1$, the values of $\ddd$ satisfy the relations in parts \ref{da} and \ref{dha}. Now we consider colours: if $d=1$, then clearly $\alpha,\hbeta$ and $\gamma$ are all grey; similarly if $d=2$ then $\halpha,\beta$ and $\hgamma$ are all grey. It remains to consider the case where $d=1$, so that $\ddd\halpha=\ddd\beta=\ddd\hgamma=0$. This means that the values $l$ and~$m$ defined above coincide. Note first that we cannot be in case \ref{low}, because then $\sigma$ would have $i$ parts in the range $\{h-i+1,\dots,h+i-1\}$, which cannot happen if $\sigma$ is an \hbc. In the other four cases, we find that
\begin{align*}
\halpha&\text{ has a $2h$-bar with leg length }2l+2\\
\beta&\text{ has two $h$-bars each with leg length }l\\
\hgamma&\text{ has a $2h$-bar with leg length }2l+1.
\end{align*}
So $\halpha,\beta,\hgamma$ are all black if $l+\bb\sigma$ is odd, and all white otherwise.

Now we prove part \ref{aghb}. By \cref{domddd}, $\alpha$, $\gamma$ and $\la$ are comparable in the dominance order, while  $\psi_i(\la)$ and $\hbeta$ are also comparable in the dominance order. Now we claim that $\la\domsby\alpha$ \iff $\psi_i(\la)\domsby\hbeta$: if $\la\domsby\alpha$, then $\la\lcolex\alpha$, so by \cref{unexlex} $\psi_i(\la)\lcolex\hbeta$, and so $\psi_i(\la)\domsby\hbeta$. On the other hand, if $\la\doms\alpha$, then $\la\glex\alpha$, so by \cref{unexlex} $\psi_i(\la)\glex\hbeta$, so $\psi_i(\la)\doms\hbeta$. So $\la\domsby\alpha$ \iff $\psi_i(\la)\domsby\hbeta$, as claimed. In exactly the same way we prove that $\la\domsby\gamma$ \iff $\psi_i(\la)\domsby\hbeta$, which gives \ref{aghb}. Part \ref{hahgb} is proved in the same way.
\end{pf}

Now we come to the special partitions in $\pb\sigma$ and their counterparts in $\pb\tau$. Our first task is to show how these partitions compare in $\pb\sigma$ and $\pb\tau$.

\begin{lemma}\label{specexc}\indent
\begin{enumerate}
\vspace{-\topsep}
\item
Suppose $\star$ is one of the symbols $\xx{}$, $\shp{}$, $\nat{}$. Then $\star_\sigma$ is unexceptional for the pair $(\blk\sigma,\blk\tau)$, and $\psi_i(\star_\sigma)=\star_\tau$.
\item
Suppose $\star$ is one of the symbols $\flt{}$, $\ppi{}$, $\yy{}$. Then one of the following happens:
\begin{enumerate}
\item
$\star_\sigma$ is unexceptional for the pair $(\blk\sigma,\blk\tau)$, and $\psi_i(\star_\sigma)=\star_\tau$.
\item
$\star_\sigma$ is exceptional for the pair $(\blk\sigma,\blk\tau)$, with $\star_\sigma=\beta$, $\star_\tau=\halpha$.
\end{enumerate}
\end{enumerate}
\end{lemma}

\begin{pf}
We consider the six special partitions separately. Recall from \cref{21stuff}\ref{21stuffexcnu} and \cref{23stuff}\ref{23stuffexcnu} that if there are exceptional partitions (i.e.\ if $i=0$ and $k=3$, or $1\ls i<n$ and $k=1$) then the exceptional partitions in $\pb\sigma$ have abacus notation $\ba{-i}$, $\ba{i,i+1}$, $\ba{i+1}$.
\begin{description}
\item[\hbox to \leftmargin{\hfil$\xx{}$ }]
As in the definition, take $1\ls a<b$ minimal such that $a,b,h-a,h-b\notin\sigma$; then $\xx\sigma=\sigma\sqcup(h-a,h-b,b,a)$, which has abacus notation $\ba{a,b}$. Given this abacus notation, $\xx\sigma$ cannot be exceptional; note that we cannot have $i=a=b-1$, because if $a=b-1$ then $\sigma$ has no addable $a$-nodes.

Now $\psi_i(\xx\sigma)$ has abacus notation $\ba{\hat a,\hat b}$, where
\[
\hat a=
\begin{cases}
a-1&\text{if $i=a-1$}\\
a+1&\text{if $i=a+1$}\\
a&\text{otherwise},
\end{cases}\qquad
\hat b=
\begin{cases}
b-1&\text{if $i=b-1$}\\
b+1&\text{if $i=b+1$}\\
b&\text{otherwise},
\end{cases}
\]
and a case-by-case check then shows that $\ba{\hat a,\hat b}=\xx\tau$.
\item[\hbox to \leftmargin{\hfil$\shp{}$ }]
In this case take $a\gs1$ minimal such that $a,h-a\notin\sigma$; then $\shp\sigma=\sigma\sqcup(h,h-a,a)$, with abacus notation $\ba{0,a}$. Hence $\shp\sigma$ is not exceptional (in the case $i=0$, observe that $a$ cannot equal $1$ because we would then have $k=1$). Now $\psi_i(\shp\sigma)=\ba{0,\hat a}$, where $\hat a$ is defined as in the previous case, and we get $\psi_i(\shp\sigma)=\shp\tau$.
\item[\hbox to \leftmargin{\hfil$\nat{}$ }]
The abacus notation for $\nat\sigma$ is $\ba{0,0}$, so $\nat\sigma$ is unexceptional. It is easily seen that $\psi_i(\nat\sigma)$ also has abacus notation $\ba{0,0}$, so equals $\nat\tau$.
\item[\hbox to \leftmargin{\hfil$\flt{}$ }]
Take $a>h$ minimal such that $a,2h-a\notin\sigma$. Then
\[
\flt\sigma=
\begin{cases}
\ba{a-h}&\text{if $a-h\notin\sigma$}\\
\ba{h-a}&\text{if $a-h\in\sigma$}.
\end{cases}
\]
Now if $i=0$, then $\flt\sigma$ cannot be exceptional (note that in this case $1\notin\sigma$, so $\flt\sigma$ cannot equal~$\ba1$). If $1\ls i<n$ and $k=1$, then $\flt\sigma$ is exceptional provided $i=a$ and $a-h\in\sigma$: then the $[2:1]$-pair is of type B, and so $\flt\sigma$ is the exceptional partition $\beta$. For every $h<c<a$ either $c$ or $2h-c$ lies in $\sigma$, so the same is true in $\tau$; in addition $a,a-h,2h-a\notin\tau$, so that $\flt\tau=\ba a$, which is the exceptional partition $\halpha$ for the pair $(\blk\sigma,\blk\tau)$.

When $\flt\sigma$ is unexceptional, it is easily checked (by a similar argument to that used in the cases above) that $\psi_i(\flt\sigma)=\flt\tau$.
\item[\hbox to \leftmargin{\hfil$\ppi{}$ }]
Take $a$ minimal such that $a\in\sigma\sqcup(h)$ but $a+h\notin\sigma$. Then
\[
\ppi\sigma=
\begin{cases}
\ba{0,a}&\text{if $a\ls n$}\\
\ba{0,h-a}&\text{if $n<a<h$}\\
\ba0&\text{if $a=h$}.
\end{cases}
\]
Clearly then $\ppi\sigma$ is unexceptional if $i\gs1$, so suppose $i=0$ and $k\gs3$. If $a\ls n$ then we would need $a=1$ in order for $\ppi\sigma$ to be exceptional; but $i=0$ implies that $1\notin\sigma$, a contradiction. If $n<a<h$ and $a=h-1$ then we get $\ppi\sigma=\ba{0,1}$, which is the exceptional partition $\beta$. Then it is easily checked that $\ppi\tau=\ba0$, which is the exceptional partition $\halpha$. Finally, if $a=h$ then the definition of $a$ means that $2h-1\in\sigma$, so that $k\gs5$, so $\ppi\sigma$ is unexceptional.

In the cases where $\ppi\sigma$ is unexceptional it is easily checked that $\psi_i(\ppi\sigma)=\ppi\tau$.
\item[\hbox to \leftmargin{\hfil$\yy{}$ }]
Take the minimal integers $a,b$ such that $a,b\notin\sigma$, $a-h\in\sigma$ and $b-h\in\sigma\sqcup(a)$, and let $\dot a$, $\dot b$ be the integers in $\{-n,\dots,-1,1,\dots,n\}$ congruent to $a$ and $b$ modulo $h$. Then
\[
\yy\sigma=
\begin{cases}
\ba{|\dot a|,|\dot b|}&\text{if $b<a+h$}\\
\ba{\dot a}&\text{if $b=a+h$}.
\end{cases}
\]
Clearly then $\yy\sigma$ is unexceptional for $i=0$. If $i\gs1$ and $b<a+h$, then $\yy\sigma$ is exceptional provided $b=a+h-1$ and $i$ equals either $\dot b$ or $-1-\dot b$ (whichever is positive). Then $(\blk\sigma,\blk\tau)$ is a $[2:1]$-pair of type A (if $\dot b>0$) or C (if $-1-\dot b>0$), and $\yy\sigma$ is the exceptional partition $\beta$, while $\yy\tau$ is obtained from $\yy\sigma$ by replacing $a-h$ with $a-h+1$, and so coincides with the exceptional partition $\halpha$.

If $b=a+h$, then $\yy\sigma$ is not exceptional: for this we would need to have $a-1\in\sigma$ but $a-1+h\notin\sigma$, but this contradicts the definition of $b$.

In the cases where $\yy\sigma$ is unexceptional it is easily checked that $\psi_i(\yy\sigma)=\yy\tau$.\qedhere
\end{description}
\end{pf}

\section{The inductive step}\label{indstep}

In this section we use the results of \cref{2ksection} to complete the proof of \cref{mainwt2}.

\begin{basicbox}
\textbf{Assumptions in force for \cref{indstep}}

$\sigma$ and $\tau$ are \hbcs, and $\blk\sigma$ and $\blk\tau$ form a $[2:k]$-pair of residue $i$, where either $i\gs1$ or  $k\gs3$.
\end{basicbox}

\begin{propn}\label{nonexc}
Suppose $\blk\sigma$ and $\blk\tau$ are Scopes--Kessar equivalent. If \cref{mainwt2} holds for $\blk\sigma$ then it holds for $\blk\tau$.
\end{propn}

\begin{pf}
By \ref{fgun} the canonical basis for $\blk\tau$ is obtained from the canonical basis for $\blk\sigma$ by applying $\psi_i$ and extending linearly. By \cref{unexddd,unexcl}, $\ddd\psi_i(\la)=\ddd\la$ and $\col\psi_i(\la)=\col\la$ for every $\la\in\pb\sigma$. Next, by \cref{unexdom}, if $\la,\mu\in\pb\sigma$ with $|\ddd\la-\ddd\mu|\ls1$ then $\la\domby\mu$ \iff $\psi_i(\la)\domby\psi_i(\mu)$. Finally, $\psi_i$ sends each special partition in $\pb\sigma$ to the corresponding special partition in $\pb\tau$, by \cref{specexc}. So if \cref{mainwt2} holds in $\blk\sigma$, then it holds in $\blk\tau$.
\end{pf}

Now we consider the more difficult case where $\blk\sigma$ and $\blk\tau$ are not Scopes--Kessar equivalent.

\begin{basicbox}
\textbf{Assumptions and notation in force for the rest of \cref{indstep}}

Either $k=1$ and $1\ls i<n$, or $k=3$ and $i=0$. $\alpha\domsby\beta\domsby\gamma$ are the exceptional partitions in $\pb\sigma$, and $\halpha\domsby\hbeta\domsby\hgamma$ the corresponding partitions in~$\blk\tau$. Let $d=\ddd\alpha=\ddd\beta+1=\ddd\gamma$.
\end{basicbox}

First we look at the partition $\mu^+$ defined in \ref{nonspecial} for a restricted \hsp $\mu\neq\xx\sigma,\shp\sigma,\nat\sigma$; this is the least dominant partition such that $\mu\domsby\mu^+$, $\ddd\mu=\ddd\mu^+$ and $\col\mu=\col\mu^+$. Part of the statement of \ref{nonspecial} is that $\mu^+$ is defined, and we address this first. For this we need to define another bijection from $\pb\sigma$ to $\pb\tau$: if $\la\in\pb\sigma$, define
\[
\isp_i(\la)=
\begin{cases}
\psi_i(\la)&\text{if $\la$ is unexceptional}\\
\hbeta&\text{if $\la=\alpha$}\\
\halpha&\text{if $\la=\beta$}\\
\hgamma&\text{if $\la=\gamma$}.
\end{cases}
\]

\begin{lemma}\label{nonspestart}
Suppose $\mu\in\blk\sigma$ is restricted, with $\mu\neq\xx\sigma,\shp\sigma,\nat\sigma$, and that $\mu^+$ is defined. Then $\psi_i(\mu)^+$ is defined and equals $\isp_i(\mu^+)$.
\end{lemma}

\begin{pf}
The case $\mu=\alpha$ is easy: here $\psi_i(\mu)=\halpha$, and from \cref{domexc} $\halpha^+$ is defined and equals $\hgamma$.

So assume that $\mu\neq\alpha$, and for the rest of the proof write $\bar\mu=\psi_i(\mu)$ and $\check\mu=\isp_i(\mu^+)$. The assumption $\mu\neq\alpha$ means that $\mu^+\neq\gamma$. Using \cref{unexddd,unexcl,unexdom,domexc} we find that $\bar\mu\domsby\check\mu$, $\ddd\bar\mu=\ddd\check\mu$ and $\col\bar\mu=\col\check\mu$. So $\bar\mu^+$ is certainly well-defined. If $\bar\mu^+\neq\check\mu$, then there is a partition $\rho\in\pb\tau$ such that $\bar\mu\domsby\rho\domsby\check\mu$, $\ddd\bar\mu=\ddd\rho$ and $\col\bar\mu=\col\rho$. We can assume $\rho\neq\hgamma$ (because if $\rho=\hgamma$, we can replace it with~$\halpha$ and it will still have the same properties), so $\rho=\isp_i(\epsilon)$ for some $\epsilon\neq\gamma$. But then we find that $\mu\domsby\epsilon\domsby\mu^+$, $\ddd\mu=\ddd\epsilon$ and $\col\mu=\col\epsilon$, contradicting the definition of~$\mu^+$.
\end{pf}

Now for a restricted partition $\mu\in\pb\sigma$ for which $\mu^+$ is defined, set
\[
A_\mu=\lset{\la\in\pb\sigma}{\mu\domsby\la\domsby\mu^+\text{ and }|\ddd\la-\ddd\mu|=1}.
\]

\begin{lemma}\label{nonspeint}
Suppose $\mu\in\blk\sigma$ is restricted, with $\mu\neq\xx\sigma,\shp\sigma,\nat\sigma$, and that $\mu^+$ is defined.
\begin{enumerate}
\item
If $\la\in\pb\sigma$ is unexceptional, then $\la\in A_\mu$ \iff $\psi_i(\la)\in A_{\psi_i(\mu)}$.
\item
$\alpha\in A_\mu\Leftrightarrow\gamma\in A_\mu\Leftrightarrow\hbeta\in A_{\psi_i(\mu)}$.
\item
$\beta\in A_\mu\Leftrightarrow\halpha\in A_{\psi_i(\mu)}\Leftrightarrow\hgamma\in A_{\psi_i(\mu)}$.
\end{enumerate}
\end{lemma}

\begin{pf}
This follows from \cref{unexddd,unexdom,domexc,nonspestart}.
\end{pf}

Now we show how to determine the canonical basis vectors for $\tau$ from those for $\si$, when $\blk\si$ and~$\blk\tau$ are not Scopes--Kessar equivalent. First we assume that \cref{mainwt2} holds for a non-special partition $\mu\in\blk\si$, and narrow down the possibilities for $\dn\al\mu,\dn\be\mu,\dn\ga\mu$.

\begin{lemma}\label{i1munon}
Suppose $i\gs1$, and $\mu$ is a non-special partition in $\blk\si$, and that \cref{mainwt2} holds for $\mu$.
\begin{enumerate}
\item\label{mnm}
$\dn\al\mu,\dn\be\mu,\dn\ga\mu\in\{0,1,q,q^2,q^3,q^4\}$.
\item\label{al1}
If $\dn\al\mu=1$, then $\dn\be\mu=q^2$ and $\dn\ga\mu=q^4$.
\item\label{ga1}
If $\dn\ga\mu=1$, then $\dn\al\mu=\dn\be\mu=0$.
\item\label{al4}
If $\dn\al\mu\in\{q^3,q^4\}$, then $\dn\be\mu=\dn\ga\mu=0$.
\item\label{sym1}
The Laurent polynomial $q^{-4}\dn\al\mu+q^{-2}\dn\be\mu+\dn\ga\mu$ is symmetric in $q$ and $q^{-1}$.
\item\label{funny1}
$(\dn\al\mu,\dn\be\mu,\dn\ga\mu)\neq(0,v,v),(v^2,v^2,v^2),(v,0,v^3),(v,v^2,v^3)$.
\end{enumerate}
\end{lemma}

\begin{pf}
(\ref{mnm}) follows from the assumption that \cref{mainwt2} holds for $\mu$. (\ref{al1}) comes from \cref{21stuff}\ref{sigma1cbv}, while (\ref{ga1},\ref{al4}) come from the statement of \cref{mainwt2} for $\mu$ and the fact that $\al\domsby\be\domsby\ga$. For (\ref{sym1}), let $\delta$ be the $h$-bar-core used in the proof of \cref{21stuff}. Then
\[
e_i\alpha=q^{-4}\delta,\qquad
e_i\beta=q^{-2}\delta,\qquad
e_i\gamma=\delta
\]
and $e_i\la=0$ for all unexceptional $\la\in\pb\sigma$. Since $\can\delta=\delta$, we obtain
\[
e_i\can\mu=\left(q^{-4}\dn\alpha\mu+q^{-2}\dn\beta\mu+\dn\gamma\mu\right)\can\delta.
\]
This vector is bar-invariant, so $q^{-4}\dn\alpha\mu+q^{-2}\dn\beta\mu+\dn\gamma\mu$ must be symmetric in $q$ and $q^{-1}$.

Finally (\ref{funny1}) comes from the statement of \cref{mainwt2} for $\mu$ together with the fact (from \cref{domexc}) that $\ddd\al=\ddd\be+1=\ddd\ga$.
\end{pf}

\begin{lemma}\label{possmunon}
Suppose $i\gs1$, and $\mu$ is a non-special partition in $\blk\si$, and that \cref{mainwt2} holds for $\mu$. Then $\dn\al\mu$, $\dn\be\mu$, $\dn\ga\mu$, $\dn{\halpha}{\psi_i(\mu)}$ $\dn{\hbeta}{\psi_i(\mu)}$, $\dn{\hgamma}{\psi_i(\mu)}$ are given by one of the rows of \cref{tablei1}.
\end{lemma}

\begin{table}[ht]
\setlength\plx{45pt}
\[
\begin{array}{|ccc|ccc|}\hline
\hbox to\plx{\hfil$\dn\alpha\mu$\hfil}&\hbox to\plx{\hfil$\dn\beta\mu$\hfil}&\hbox to\plx{\hfil$\dn\gamma\mu$\hfil}&\hbox to\plx{\hfil$\dn{\halpha}{\psi_i(\mu)}$\hfil}&\hbox to\plx{\hfil$\dn{\hbeta}{\psi_i(\mu)}$\hfil}&\hbox to\plx{\hfil$\dn{\hgamma}{\psi_i(\mu)}$\hfil}\\\hline
0&0&0&0&0&0\\
0&0&1&0&1&q^2\\
0&1&q^2&0&0&1\\
0&q^2&0&q^2&0&q^2\\
1&q^2&q^4&1&q^2&q^4\\
q^2&0&q^2&0&q^2&0\\
q^2&q^4&0&q^4&0&0\\
q^4&0&0&q^2&q^4&0\\
\hline
\end{array}
\]
\caption{}\label{tablei1}
\end{table}

\begin{pf}
It is a simple matter to check that the only possible triples $(\dn\al\mu,\dn\be\mu,\dn\ga\mu)$ satisfying all the conditions of \cref{i1munon} are those given in the table. In each case, we can compute the bar-invariant vector $f_i\can\mu$ using \ref{fun} and part~\ref{21stuffpsi} of the present \lcnamecref{21stuff}, and then reduce this using the known vector $\can\halpha$ to obtain a \cbv in $\blk\tau$, which turns out to be $\can{\psi_i(\mu)}$. We give an example of this calculation.

Suppose we are in the third case in \cref{tablei1}, with
\[
\can\mu=\beta+q^2\gamma+\sum_{\la\text{ unexceptional}}\dn\la\mu\la,
\]
so that in particular $\mu=\beta$. Then
\[
f_i\can\mu=\halpha+q^2\hbeta+(q^4+1)\hgamma+\sum_{\la\text{ unexceptional}}\dn\la\mu\psi_i(\la).
\]
Subtracting $\can\halpha$, we obtain
\[
\hgamma+\sum_{\la\text{ unexceptional}}\dn\la\mu\psi_i(\la)
\]
which therefore equals $\can\hgamma=\can{\psi_i(\mu)}$.
\end{pf}

Now we do the same for $i=0$. The next two lemmas are proved in the same way as \cref{i1munon,possmunon}. For part (\ref{sym0}) of \cref{i0munon}, we use
\[
e_0\alpha=q^{-4}\delta,\qquad
e_0\beta=(q^{-1}+q^{-3})\delta,\qquad
e_0\gamma=(q+q^{-1})\delta
\]
where $\delta=\si\sqcup(2h-1)$.

\begin{lemma}\label{i0munon}
Suppose $i=0$ and $\mu$ is a non-special partition in $\blk\si$, and that \cref{mainwt2} holds for $\mu$.
\begin{enumerate}
\item
$\dn\al\mu,\dn\be\mu,\dn\ga\mu\in\{0,1,q,q^2,q^3,q^4\}$.
\item
If $\dn\ga\mu=1$, then $\dn\al\mu=\dn\be\mu=0$.
\item
If $\dn\al\mu\in\{q^3,q^4\}$, then $\dn\be\mu=\dn\ga\mu=0$.
\item\label{sym0}
The Laurent polynomial $q^{-4}\dn\al\mu+(q^{-1}+q^{-3})\dn\be\mu+(q+q^{-1})\dn\ga\mu$ is symmetric in $q$ and $q^{-1}$.
\item
$(\dn\al\mu,\dn\be\mu,\dn\ga\mu)\neq(0,v,v),(v^2,v^2,v)$.
\end{enumerate}
\end{lemma}

\begin{lemma}\label{0possmunon}
Suppose $i=0$ and $\mu$ is a non-special partition in $\blk\si$, and that \cref{mainwt2} holds for $\mu$. Then $\dn\al\mu$, $\dn\be\mu$, $\dn\ga\mu$, $\dn{\halpha}{\psi_i(\mu)}$ $\dn{\hbeta}{\psi_i(\mu)}$, $\dn{\hgamma}{\psi_i(\mu)}$ are given by one of the rows of \cref{tablei0}.
\end{lemma}

\begin{table}[ht]
\setlength\plx{45pt}
\[
\begin{array}{|ccc|ccc|}\hline
\hbox to\plx{\hfil$\dn\alpha\mu$\hfil}&\hbox to\plx{\hfil$\dn\beta\mu$\hfil}&\hbox to\plx{\hfil$\dn\gamma\mu$\hfil}&\hbox to\plx{\hfil$\dn{\halpha}{\psi_i(\mu)}$\hfil}&\hbox to\plx{\hfil$\dn{\hbeta}{\psi_i(\mu)}$\hfil}&\hbox to\plx{\hfil$\dn{\hgamma}{\psi_i(\mu)}$\hfil}\\\hline
0&0&0&0&0&0\\
0&0&1&0&1&q^2\\
0&1&q^2&0&0&1\\
0&q^2&0&q^2&0&q^2\\
1&q&q^3&1&q^2&q^4\\
q^2&0&q&0&q&0\\
q^2&q^3&0&q^3&0&0\\
q^4&0&0&q&q^3&0\\
\hline
\end{array}
\]
\caption{}\label{tablei0}
\end{table}

Now we can complete the inductive step for non-special partitions.

\begin{propn}\label{nonspe1}
Suppose $\mu\in\blk\sigma$ is restricted, with $\mu\neq\xx\sigma,\shp\sigma,\nat\sigma$, and that \ref{nonspecial} holds for $\mu$. Then \ref{nonspecial} holds for $\psi_i(\mu)$.
\end{propn}

\begin{pf}
First suppose $\mu=\alpha$. Then $\psi_i(\mu)=\halpha$, and we know $\can\halpha$ from \cref{21stuff}\ref{sigma1cbv} or \cref{23stuff}\ref{gminusi}. From \cref{domexc} $\can\halpha$ satisfies \cref{mainwt2}.

Now take $\mu\neq\alpha$. Then $\dn\alpha\mu$, $\dn\beta\mu$, $\dn\gamma\mu$, $\dn\alpha\mu$, $\dn\beta\mu$, $\dn\gamma\mu$ are given by one of the rows of \cref{tablei1} (if $i\gs1$) or \cref{tablei0} (if $i=0$). Now the result follows from \cref{21stuff}\ref{21stuffpsi} or \cref{21stuff}\ref{21stuffpsi}, together with \cref{nonspestart} and \cref{nonspeint}.
\end{pf}

\begin{eg}
We give an example to illustrate the proof of \cref{nonspe1}. Suppose $\mu=\gamma$, so that $(\dn\alpha\gamma,\dn\beta\gamma,\dn\gamma\gamma)=(0,0,1)$. Then $\psi_i\mu=\hbeta$, and by \cref{possmunon} or \cref{0possmunon} $(\dn\halpha\hbeta,\dn\hbeta\hbeta,\dn\hgamma\hbeta)=(0,1,q^2)$ while $\dn{\psi_i(\la)}\hbeta=\dn\la\gamma$ for unexceptional $\la$. Now $\gamma^+$ is unexceptional, so $\hat\beta^+=\psi_i(\gamma^+)$ by \cref{nonspestart}. Now by \cref{unexddd,unexdom,domexc}
\[
A_{\hat\beta}=\lset{\psi_i(\la)}{\la\in A_\gamma}\cup\{\gamma\}.
\]
Hence the coefficients in $\can\hbeta$ satisfy \cref{mainwt2}.
\end{eg}

Now we come to the special partitions $\xx\sigma$, $\shp\sigma$ and $\nat\sigma$. Because $\bb\sigma=\bb\tau$, the special partitions that are defined in $\pb\sigma$ are also defined in $\pb\tau$ and vice versa; moreover, $\nat\sigma$ is restricted \iff $\nat\tau$ is restricted.

\begin{propn}\label{indxx}
Suppose $\bb\sigma\ls n-2$ and \cref{mainwt2} holds for $\mu=\xx\sigma$. Then it holds for $\xx\tau$.
\end{propn}

\begin{pf}
We know from \cref{specexc} that $\xx\sigma$, $\shp\sigma$ and $\nat\sigma$ are unexceptional, with $\psi_i(\xx\sigma)=\xx\tau$, $\psi_i(\shp\sigma)=\shp\tau$, $\psi_i(\nat\sigma)=\nat\tau$. In addition, either $\flt\sigma$ is unexceptional and $\flt\tau=\psi_i(\flt\sigma)$, or $\flt\sigma=\beta$ and $\flt\tau=\halpha$.

Now let
\begin{align*}
D(\sigma)&=\lset{\la}{\xx\sigma\domsby\la\domsby\shp\sigma,\ \ddd\la=2}\\
E(\sigma)&=\lset{\la}{\shp\sigma\domsby\la\domsby\flt\sigma,\ \ddd\la=2}\\
F(\sigma)&=D(\sigma)\cup E(\sigma)\cup\{\xx\si,\shp\si,\nat\sigma,\flt\sigma\},
\end{align*}
and define $D(\tau)$, $E(\tau)$, $F(\tau)$ similarly. By considering $\ddd$-values and the dominance order, we find six possibilities for the intersection of $\{\alpha,\beta,\gamma\}$ with $F(\sigma)$.
\begin{enumerate}
\item
$\alpha,\beta,\gamma\notin F(\sigma)$. In this case $D(\tau)=\psi_i(D(\sigma))$ and $E(\tau)=\psi_i(E(\sigma))$.
\item
$\beta\in D(\sigma)$, $\alpha,\gamma\notin F(\sigma)$.
\item
$\alpha,\gamma\in D(\sigma)$, $\beta\notin F(\sigma)$.
\item
$\beta\in E(\sigma)$, $\alpha,\gamma\notin F(\sigma)$. In this case $D(\tau)=\psi_i(D(\sigma))$ and $E(\tau)=\psi_i(E(\sigma))\cup\{\halpha\}$.
\item\label{xxbe}
$\alpha,\gamma\in E(\sigma)$, $\beta\notin F(\sigma)$. In this case $D(\tau)=\psi_i(D(\sigma))$ and $E(\tau)=\psi_i(E(\sigma))\setminus\{\halpha\}$.
\item\label{xxae}
$\alpha\in E(\sigma)$, $\beta=\flt\sigma$, $\gamma\notin F(\sigma)$. In this case $D(\tau)=\psi_i(D(\sigma))$, $E(\tau)=\psi_i(E(\sigma))\setminus\{\halpha\}$ and $\flt\tau=\halpha$.
\end{enumerate}
In fact cases (2) and (3) cannot occur. This follows from the same argument used to prove \cref{i1munon}(6): in case (2) we would have $(\dn\al{\xx\si},\dn\be{\xx\si},\dn\ga{\xx\si})=(0,q,0)$, and applying $e_i$ to $G(\xx\si)$ would give a non-bar-invariant vector; similarly for case (3). By the same reasoning, cases (5) and (6) can only happen when $i\gs1$.

In cases (1,4,5,6), the values $\dn\la{\xx\tau}$ can be calculated using the technique in the proof of \cref{possmunon}: applying $f_i^{(k)}$ to $G(\xx\si)$, and subtracting a suitable multiple of $G(\halpha)$ to obtain $G(\xx\tau)$. In each case we find that \cref{mainwt2} holds for $\xx\tau$.
\end{pf}

\begin{propn}\label{indshp}
Suppose $\bb\sigma\ls n-1$ and \cref{mainwt2} holds for $\mu=\shp\sigma$. Then it holds for $\shp\tau$.
\end{propn}

\begin{pf}
The structure of the proof is exactly as for \cref{indxx}. We know from \cref{specexc} that $\shp\sigma$ and $\nat\sigma$ are unexceptional, with $\psi_i(\shp\sigma)=\shp\tau$, $\psi_i(\nat\sigma)=\nat\tau$. In addition, either $\flt\sigma$ is unexceptional with $\flt\tau=\psi_i(\flt\si)$ or $\flt\si=\be$ and $\flt\tau=\halpha$; a similar statement holds for $\ppi\sigma$ and $\ppi\tau$.

Now let
\begin{align*}
D(\sigma)&=\lset{\la}{\shp\sigma\domsby\la\domsby\flt\sigma,\ \ddd\la=2}\\
E(\sigma)&=\lset{\la}{\flt\sigma\domsby\la\domsby\ppi\sigma,\ \ddd\la=1}\\
F(\sigma)&=D(\sigma)\cup E(\sigma)\cup\{\shp\si,\nat\si,\flt\sigma,\ppi\sigma\},
\end{align*}
and define $D(\tau)$, $E(\tau)$, $F(\tau)$ similarly. Now we find seven possibilities for the intersection of $\{\alpha,\beta,\gamma\}$ with $F(\sigma)$.
\begin{enumerate}
\item
$\alpha,\beta,\gamma\notin F(\sigma)$. In this case $D(\tau)=\psi_i(D(\sigma))$ and $E(\tau)=\psi_i(E(\sigma))$.
\item
$\beta\in D(\sigma)$, $\alpha,\gamma\notin F(\sigma)$. In this case $D(\tau)=\psi_i(D(\sigma))\cup\{\halpha\}$ and $E(\tau)=\psi_i(E(\sigma))$.
\item\label{shagd}
$\alpha,\gamma\in D(\sigma)$, $\beta\notin F(\sigma)$. In this case $D(\tau)=\psi_i(D(\sigma))\setminus\{\halpha\}$, and $E(\tau)=\psi_i(E(\sigma))$.
\item\label{shad}
$\alpha\in D(\sigma)$, $\beta=\flt\sigma$, $\gamma\notin F(\sigma)$. In this case $D(\tau)=\psi_i(D(\sigma))\setminus\{\halpha\}$, $\flt\tau=\halpha$ and $E(\tau)=\psi_i(E(\sigma))\cup\{\hgamma\}$.
\item
$\beta\in E(\sigma)$, $\alpha,\gamma\notin F(\sigma)$. In this case $D(\tau)=\psi_i(D(\sigma))$ and $E(\tau)=\psi_i(E(\sigma))\cup\{\halpha\}$.
\item\label{shage}
$\alpha,\gamma\in E(\sigma)$, $\beta\notin F(\sigma)$. In this case $D(\tau)=\psi_i(D(\sigma))$ and $E(\tau)=\psi_i(E(\sigma))\setminus\{\halpha\}$.
\item\label{shae}
$\alpha\in E(\sigma)$, $\beta=\ppi\sigma$, $\gamma\notin F(\sigma)$. In this case $D(\tau)=\psi_i(D(\sigma))$, $E(\tau)=\psi_i(E(\sigma))\setminus\{\halpha\}$ and $\ppi\tau=\halpha$.
\end{enumerate}
Applying $e_i$, we find that cases \ref{shagd}, \ref{shad} and \ref{shage} can only happen when $i\gs1$, while case \ref{shae} can only happen if $i=0$. Now as in the proof of \cref{indxx} we can find $G(\shp\tau)$ and verify that \cref{mainwt2} holds for $\shp\tau$.
\end{pf}

\begin{propn}\label{indnat}
Suppose $\bb\sigma\gs1$ and \cref{mainwt2} holds for $\mu=\nat\sigma$. Then it holds for $\nat\tau$.
\end{propn}

\begin{pf}
We use the same technique as in \cref{indxx,indshp}. This time we define
\begin{align*}
D(\sigma)&=\lset{\la}{\nat\sigma\domsby\la\domsby\ppi\sigma,\ \ddd\la=1}\\
E(\sigma)&=\lset{\la}{\ppi\sigma\domsby\la\domsby\yy\sigma,\ \ddd\la=1}\\
F(\sigma)&=D(\sigma)\cup E(\sigma)\cup\{\ppi\sigma,\yy\sigma\},
\end{align*}
and we find seven possibilities for the intersection of $\{\alpha,\beta,\gamma\}$ with $F(\sigma)$.
\begin{enumerate}
\item
$\alpha,\beta,\gamma\notin F(\sigma)$. In this case $D(\tau)=\psi_i(D(\sigma))$ and $E(\tau)=\psi_i(E(\sigma))$.
\item
$\beta\in D(\sigma)$, $\alpha,\gamma\notin F(\sigma)$. In this case $D(\tau)=\psi_i(D(\sigma))\cup\{\halpha\}$ and $E(\tau)=\psi_i(E(\sigma))$.
\item\label{ntagd}
$\alpha,\gamma\in D(\sigma)$, $\beta\notin F(\sigma)$. In this case $D(\tau)=\psi_i(D(\sigma))\setminus\{\halpha\}$, and $E(\tau)=\psi_i(E(\sigma))$.
\item\label{ntadge}
$\alpha\in D(\sigma)$, $\beta=\ppi\sigma$, $\gamma\in E(\sigma)$. In this case $D(\tau)=\psi_i(D(\sigma))\setminus\{\halpha\}$, $\ppi\tau=\halpha$ and $E(\tau)=\psi_i(E(\sigma))$.
\item
$\beta\in E(\sigma)$, $\alpha,\gamma\notin F(\sigma)$. In this case $D(\tau)=\psi_i(D(\sigma))$ and $E(\tau)=\psi_i(E(\sigma))\cup\{\halpha\}$.
\item\label{ntage}
$\alpha,\gamma\in E(\sigma)$, $\beta\notin F(\sigma)$. In this case $D(\tau)=\psi_i(D(\sigma))$ and $E(\tau)=\psi_i(E(\sigma))\setminus\{\halpha\}$.
\item\label{ntae}
$\alpha\in E(\sigma)$, $\beta=\yy\sigma$, $\gamma\notin F(\sigma)$. In this case $D(\tau)=\psi_i(D(\sigma))$, $E(\tau)=\psi_i(E(\sigma))\setminus\{\halpha\}$ and $\yy\tau=\halpha$.
\end{enumerate}
(Cases \ref{ntagd}, \ref{ntage} and \ref{ntae} can only happen when $i\gs1$, while case \ref{ntadge} can only happen if $i=0$.) Now we proceed as in \cref{indxx}.
\end{pf}

Finally we can complete the proof of the main theorem.

\begin{pf}[Proof of \cref{mainwt2}]
We proceed by induction on $\card\tau$. Let $l=\len\tau$, and consider the three possibilities in \cref{non2core}. If $\tau=(l,l-1,\dots,1)$ with $l\ls n$, then the results of \cref{basesec} show that \cref{mainwt2} holds for $\tau$. Alternatively, there is a residue $i$ such that either $i\neq0$ and $\tau$ has a removable $i$-node, or $i=0$ and $\tau$ has at least three removable $i$-nodes. So define the \hbc $\sigma$ by removing all the removable $i$-nodes from $\tau$. Now the result follows from \cref{nonexc,nonspe1,indxx,indshp,indnat}, and the inductive hypothesis.
\end{pf}

\section{The Fock space of type $A^{(2)}_h$}\label{otherfocksec}

The subject of this paper is the determination of \cbcs for the $q$-deformed Fock space in type $A^{(2)}_{h-1}$, where $h\gs3$ is odd. In this section we briefly discuss the corresponding problem in type $A^{(2)}_h$. These two Kac--Moody types are related by folding of Dynkin diagrams, as illustrated below.
\[
\begin{array}{c@{\qquad}c@{\qquad}c}
\text{Type $A^{(2)}_h$}&
\begin{tikzpicture}[scale=.6,baseline=0]
\foreach\y in{-.1,.1}\draw(10,\y)--++(2,0);
\foreach\x in{2,8}\draw(\x,0)--++(2,0);
\draw(2,0)++(135:2)--++(315:2)--++(225:2);
\draw[thick,loosely dotted](5,0)--++(2,0);
\foreach\x in {2,4,8,10,12}\draw[black,fill=white](\x,0)circle(.175);
\draw[black,fill=white](2,0)++(135:2)circle(.175);
\draw[black,fill=white](2,0)++(225:2)circle(.175);
\foreach\x in{11.15}\draw(\x,-.3)--++(-.3,.3)--++(.3,.3);
\foreach\x in {1,2}\draw(2*\x,-.75)node{\small$\x$};
\draw(2,-.75)++(135:2)node{\small$\hat0$};
\draw(2,-.75)++(225:2)node{\small$\check0$};
\draw(8,-.75)node{$n-2$}++(2,0)node{$n-1$}++(2,0)node{$n$};
\end{tikzpicture}
&
\begin{tikzpicture}[scale=.6,baseline=0]
\draw(2,0)++(45:.1)--++(135:2);
\draw(2,0)++(45:-.1)--++(135:2);
\draw(2,0)++(135:.1)--++(225:2);
\draw(2,0)++(135:-.16)--++(225:2);
\draw(2,0)++(135:1.15)--++(.424,0);
\draw(2,0)++(135:1.15)--++(0,-.424);
\draw(2,0)++(225:1.15)--++(.424,0);
\draw(2,0)++(225:1.15)--++(0,.424);
\draw[black,fill=white](2,0)circle(.175);
\draw[black,fill=white](2,0)++(135:2)circle(.175);
\draw[black,fill=white](2,0)++(225:2)circle(.175);
\draw(2,-.75)node{\small$1$};
\draw(2,-.75)++(135:2)node{\small$\hat0$};
\draw(2,-.75)++(225:2)node{\small$\check0$};
\end{tikzpicture}
\\[52pt]
\text{Type $A^{(2)}_{h-1}$}&
\begin{tikzpicture}[scale=.6,baseline=0]
\foreach\y in{-.1,.1}\foreach\x in {0,10}\draw(\x,\y)--++(2,0);
\foreach\x in{2,8}\draw(\x,0)--++(2,0);
\draw[thick,loosely dotted](5,0)--++(2,0);
\foreach\x in {0,2,4,8,10,12}\draw[black,fill=white](\x,0)circle(.175);
\foreach\x in{1.15,11.15}\draw(\x,-.3)--++(-.3,.3)--++(.3,.3);
\foreach\x in {0,1,2}\draw(2*\x,-.75)node{\small$\x$};
\draw(8,-.75)node{$n-2$}++(2,0)node{$n-1$}++(2,0)node{$n$};
\end{tikzpicture}
&
\begin{tikzpicture}[scale=.6,baseline=0]
\foreach\y in{-.05,.05,-.15,.15}\draw(0,\y)--++(2,0);
\foreach\x in {0,2}\draw[black,fill=white](\x,0)circle(.175);
\foreach\x in{1.15}\draw(\x,-.3)--++(-.3,.3)--++(.3,.3);
\foreach\x in {0,1}\draw(2*\x,-.75)node{\small$\x$};
\end{tikzpicture}
\\[28pt]
&
h\gs5
&
h=3
\end{array}
\]
(In the Kac--Moody classification, type $A^{(2)}_3$ is usually referred to as $D^{(2)}_3$.)

Following the general construction for all classical types in \cite{kk} in terms of Young walls, the Fock space can be described combinatorially; we give a brief summary.  Define a \emph{decorated \hsp} to be an \hsp in which each non-zero part divisible by $h$ is decorated with an accent $\hat{}$ or~$\check{}$, and two consecutive equal parts must have opposite decorations. The $q$-deformed Fock space has the set of all decorated $h$-strict partitions as a basis, and the actions of the generators $e_i$ and $f_i$ can be described in terms of adding and removing $i$-nodes (with a suitable definition of $i$-node for $i=\hat0$ or~$\check0$). The weight spaces can be defined in term of ``decorated \hbcs'', leading to a suitable notion of \baw.

It appears that analogues of \cref{mainwt1,mainwt2} hold in this setting, and can be proved by the same techniques; in fact, the situation for weight spaces of \baw $2$ in type $A^{(2)}_h$ is simpler in several ways:
\begin{itemize}
\item
there is no exceptional behaviour for the special partitions, so a more direct analogue of Richards's theorem is possible;
\item
the \cbcs are all equal to $0$, $1$, $q$ or $q^2$, so there are fewer cases to consider in the analysis of $[2:k]$-pairs;
\item
$[2:1]$-pairs of residue $\hat0$ or $\check0$ are much more tractable, so that only the case $\tau=\varnothing$ is needed as a base case.
\end{itemize}

The relationship between the quantum groups of types $A^{(2)}_h$ and $A^{(2)}_{h-1}$ means that the \cbcs in the two types are very closely related. The ``folding'' process involved in the transition introduces the exceptional behaviour for the special partitions. We illustrate this in \cref{a4a5} by giving an example for $h=5$, showing the \cbcs for the weight spaces corresponding to the $5$-bar-core $(1)$, with the special partitions labelled.
\begin{figure}[ht]
\[
\begin{array}{c@{\qquad}c}
\begin{array}{c|cccccc|}
&\rt{5,3,2,1}&\rt{\hat5,\check5,1}&\rt{\check5,\hat5,1}&\rt{6,3,2}&\rt{6,4,1}&\rt{7,3,1}\\
\hline
(5,3,2,1)&1&\cdot&\cdot&\cdot&\cdot&\cdot\\
(\hat5,\check5,1)&\cdot&1&\cdot&\cdot&\cdot&\cdot\\
(\check5,\hat5,1)&q&\cdot&1&\cdot&\cdot&\cdot\\
(6,3,2)&q&\cdot&\cdot&1&\cdot&\cdot\\
(6,4,1)&q^2&q&q&q&1&\cdot\\
(6,5)&\cdot&q^2&\cdot&\cdot&q&\cdot\\
(7,3,1)&\cdot&\cdot&\cdot&q^2&q&1\\
(8,2,1)&\cdot&\cdot&\cdot&\cdot&q&q^2\\
(10,1)&\cdot&\cdot&q&\cdot&q^2&\cdot\\
(11)&\cdot&\cdot&q^2&\cdot&\cdot&\cdot\\\hline
\end{array}
&
\begin{array}{rc|ccccc|}
&&\rt{5,3,2,1}&\rt{5^2,1}&\rt{6,3,2}&\rt{6,4,1}&\rt{7,3,1}\\
\hline
\shp{}&(5,3,2,1)&1&\cdot&\cdot&\cdot&\cdot\\[\jmp]
\nat{}&(5^2,1)&q&1&\cdot&\cdot&\cdot\\[\jmp]
&(6,3,2)&q^2&\cdot&1&\cdot&\cdot\\
\flt{}&(6,4,1)&q^4{+}q^2&q^3{+}q&q^2&1&\cdot\\
\ppi{}&(6,5)&q^3&q^2&\cdot&q&\cdot\\
&(7,3,1)&\cdot&\cdot&q^4&q^2&1\\
&(8,2,1)&\cdot&\cdot&\cdot&q^2&q^4\\
&(10,1)&\cdot&q^2&\cdot&q^3&\cdot\\
\yy{}&(11)&\cdot&q^4&\cdot&\cdot&\cdot\\\hline
\end{array}
\\\\[0pt]
\text{type $A^{(2)}_5$}
&
\text{type $A^{(2)}_4$}
\end{array}
\]
\caption{Canonical basis coefficients in types $A_5^{(2)}$ and $A_4^{(2)}$}\label{a4a5}
\end{figure}
We can see that the first two columns of the right-hand matrix (giving the \cbvs $G(\shp{})$ and $G(\nat{})$, for which the exceptional behaviour occurs) are obtained by adding the first three columns of the left-hand matrix in pairs, with an (as yet mysterious) adjustment to the powers of $q$ occurring.

It therefore appears that a promising approach to finding \cbcs in type $A^{(2)}_{h-1}$ (and to the spin decomposition number problem for symmetric groups in characteristic $h$) is to work in type $A^{(2)}_h$ first, and understand how folding affects decomposition numbers.

\section{Application to spin representations of symmetric groups}\label{spinsec}

This paper is motivated by the decomposition number problem for spin representations of symmetric groups. Here we briefly summarise the background, and discuss the implications of our results.

Take $m\gs4$, and let $\hs m$ denote one of the two Schur covers of $\sss m$. Any representation of $\sss m$ lifts to a representation of $\hs m$; the irreducible representations which do not come from $\sss m$ in this way are called \emph{spin} representations of $\sss m$. Given a representation (or character), the \emph{associate} representation is obtained by tensoring with the (lift of the) sign representation of $\sss m$.

The classification of irreducible characters of $\hs m$ over $\bbc$ goes back to Schur \cite{schu} (though construction of the actual representations was achieved much later, by Nazarov \cite{naz}), and can be stated as follows. Say that a strict partition $\la$ is \emph{even} or \emph{odd} as the number of positive even parts of $\la$ is even or odd. For each even strict partition $\la$ of $m$, there is an irreducible self-associate character $\chi(\la)$ of $\hs m$ if $\la$ is even, and a pair of associate irreducible characters $\chi(\la)_+,\chi(\la)_-$ if $\la$ is odd. These characters are pairwise distinct, and yield all the ordinary irreducible spin characters of $\hs m$.

The classification of irreducible modular representations is due to Brundan and Kleshchev \cite{bk1,bk2}. Suppose $h$ is an odd prime, and let $\bbf$ be a field of characteristic $h$ which is a splitting field for~$\hs m$. Say that a partition is \emph{$h$-even} or \emph{$h$-odd} as the number of nodes of non-zero residue is even or odd. Then for each restricted $h$-strict partition $\mu$ of $m$, there is a self-associate irreducible Brauer character~$\phi(\mu)$ if $\mu$ is $h$-even, and a pair of associate irreducible Brauer characters $\phi(\mu)_+,\phi(\mu)_-$ if $\mu$ is $h$-odd. These Brauer characters are distinct, and give all the irreducible $h$-modular spin Brauer characters of~$\hs m$.

Two characters (ordinary or modular) lie in the same $h$-block of $\hs m$ \iff the labelling partitions have the same \hbc (except in the case of an odd partition $\la$ of \hbw $0$, where~$\chi(\la)_+$ and $\chi(\la)_-$ lie in separate simple blocks). So (apart from this slight caveat for \hbw $0$) $h$-blocks correspond precisely to the blocks (i.e.\ weight spaces in $\scrf$) studied in this paper. The defect group of a block with \baw $w$ is isomorphic to a Sylow $h$-subgroup of $\sss{wh}$; in particular, for blocks with \baw less than $h$, the defect group is abelian and the defect coincides with the \baw.

The \emph{spin decomposition number problem} then asks for the decomposition of $\chi(\la)$ or $\chi(\la)_\pm$ as a sum of irreducible Brauer characters $\phi(\mu)$ or $\phi(\mu)_\pm$. A close approximation to this problem is to consider the \emph{reduced decomposition number} obtained by combining the indecomposable projective characters corresponding to associate Brauer characters. We use a slightly unusual convention for this: for a strict partition $\la$ of $m$ and a restricted $h$-strict partition $\mu$ of $m$ we define the reduced decomposition number
\[
\rdn\la\mu=
\begin{cases}
\bp{\chi(\la)}{\phi(\mu)}
&\text{if $\la$ is even and $\mu$ is $h$-even}\\
\sqrt2\bp{\chi(\la)}{\phi(\mu)_+}
&\text{if $\la$ is even and $\mu$ is $h$-odd}\\
\sqrt2\bp{\chi(\la)_+}{\phi(\mu)}
&\text{if $\la$ is odd and $\mu$ is $h$-even}\\
\bp{\chi(\la)_+}{\phi(\mu)_+}+\bp{\chi(\la)_+}{\phi(\mu)_-}
&\text{if $\la$ is odd and $\mu$ is $h$-odd}.
\end{cases}
\]
A conjecture due to Leclerc and Thibon \cite[Conjecture~6.2]{lt} says that if $h$ is large relative to $m$ then~$\rdn\la\mu$ is determined by the integer obtained by evaluating $\dn\la\mu$ at $q=1$; specifically, if we define $n_h(\la)$ to be the number of positive parts of $\la$ divisible by $h$, then $\rdn\la\mu=2^{n_h(\la)/2}\dn\la\mu(1)$.

In fact, the original Leclerc--Thibon conjecture asserts that this relationship should hold whenever $m<h^2$. A reasonable extension to a blockwise version would say that the formula should hold in all blocks of \baw less than $h$ (regardless of $m$), i.e.\ for all blocks with abelian defect group; this is analogous to the blockwise form of James's conjecture for decomposition numbers of symmetric groups. \cref{mainwt1} and M\"uller's work \cite{muller} show that the Leclerc--Thibon conjecture is true for blocks of \baw $1$. However, the original version of the Leclerc--Thibon conjecture cannot be true: as pointed out by Tsuchioka (in private communication) it predicts that $\rdn{(11,4,3,2)}{(5^2,4,3,2,1)}=-1$ in characteristic $5$, which is absurd. Nevertheless, it seems likely that the blockwise version does hold for blocks of \baw $2$; the decomposition numbers are known for all $m\ls 18$ thanks to Maas \cite{maas}, and the conjecture can be checked in these cases. So it appears that setting $q=1$ in \cref{mainwt2} gives a formula for the reduced decomposition numbers for defect $2$ spin blocks of symmetric groups. We hope to prove this in future work.

\end{document}